\documentclass[12pt]{article}
\usepackage{fancyhdr}
\usepackage{graphicx}
\usepackage{amscd, amssymb}
\usepackage{enumerate}
\usepackage{hyperref}
\usepackage{lscape}
\newenvironment{eq}{\begin{equation}}{\end{equation}}
\newenvironment{myproof}{{\bf Proof}:}{\vskip 5mm }

\newtheorem{proposition}{Proposition}[subsection]
\newtheorem{lemma}[proposition]{Lemma}
\newtheorem{definition}[proposition]{Definition}
\newtheorem{theorem}[proposition]{Theorem}
\newtheorem{cor}[proposition]{Corollary}

\newtheorem{remark}[proposition]{Remark}

\newtheorem{cond}[proposition]{Conditions}

\newtheorem{problem}[proposition]{Problem}
\newtheorem{construction}[proposition]{Construction}
\newcommand{\llabel}[1]{\label{#1}}
\newcommand{\comment}[1]{}
\newcommand{\sr}{\rightarrow}
\newcommand{\lr}{\longrightarrow}

\newcommand{\nn}{{\bf N\rm}}

\newcommand{\uu}{\underline}

\newcommand{\wt}{\widetilde}

\newcommand{\BB}{{\bullet}}

\newcommand{\dd}{\diamond}

{\vskip
3mm}

\newcommand{\spc}{{\,\,\,\,\,\,\,}}

\comment{\pagestyle{fancy}

\renewcommand{\sectionmark}[1]{}
\renewcommand{\chaptermark}[1]{}}

\begin{document}
\parskip = 2mm
\begin{center}
{\bf\Large Martin-Lof identity types in the C-systems defined by\\
\vskip 2mm
 a universe category\footnote{\em 2000 Mathematical Subject Classification: 
18D99, 
03B15, 
18D15 
}}

\vspace{3mm}

{\large\bf Vladimir Voevodsky}\footnote{School of Mathematics, Institute for Advanced Study,
Princeton NJ, USA. e-mail: vladimir@ias.edu}$^,$\footnote{Work on this paper was supported by NSF grant 1100938.}
\vspace {3mm}

{May 2015}  
\end{center}
\begin{abstract}
This paper continues the series of papers that develop a new approach to syntax and semantics of dependent type theories. Here we study the interpretation of the rules of the identity types in the intensional Martin-Lof type theories on the C-systems that arise from universe categories.  In the first part of the paper we develop constructions that produce interpretations of these rules from certain structures on universe categories while in the second we study the functoriality of these constructions with respect to functors of universe categories. The results of the first part of the paper play a crucial role in the construction of the univalent model of type theory in simplicial sets. 
\end{abstract}

\vskip 4mm
\tableofcontents

\comment{\begin{minipage}{60mm}
He that delivereth knowledge desireth to deliver it in such form as may be soonest believed and not as may be easiest examined.

``On the Impediments of Knowledge'', from Valerius Terminus by Francis Bacon. 
\end{minipage}}

\section{Introduction}

The concept of a C-system in its present form was introduced in \cite{Csubsystems}. The type of the C-systems is constructively equivalent to the type of contextual categories defined by Cartmell in \cite{Cartmell1} and \cite{Cartmell0} but the definition of a C-system is slightly different from the Cartmell's foundational definition.

In the past decade or more, it has been a tradition in the study of type theories to consider, as the main mathematical object associated with a type theory, not a C-system by a category with families (see \cite{Dybjer}). As was observed recently all of the constructions of \cite{Cfromauniverse}, \cite{fromunivwithPi}  and of the present paper (but not of \cite{Csubsystems} or \cite{Cofamodule}!) can be either used directly or reformulated in a very straightforward way to provide very similar results for categories with families. This modification will be discussed in a separate paper or papers. 

In this introductory explanation we will distinguish between the syntactic and semantic C-systems. By a syntactic C-system we will mean a C-system that is a regular sub-quotient of a C-system of the form $CC(R,LM)$ where $R$ is a monad on sets and $LM$ is a left module over $R$, see \cite{Cofamodule} and \cite{Csubsystems}. In particular, the C-systems of all of the various versions of dependent type theory of Martin-Lof ``genus" are syntactic type systems in the sense of this definition. 

By a semantic C-system we will mean a C-system whose underlying category is a full subcategory in a category of ``mathematical'' nature such as the category of sets or the category of sheaves of sets. 

Usually one knows some good properties (i.e. consistency) of a given semantic C-system and tries to prove similar good properties of a syntactic C-system by constructing a functor from the syntactic one to the semantic one. 

To construct such a functor one tries to show that the syntactic C-system is an initial one among C-systems equipped with some collection of additional operations and then to construct operations of the required form on the semantic one. A pioneering example of this approach can be found in \cite{Streicher}.

In this paper we investigate the set of three interconnected operations on C-systems that, in the case of the syntactic C-systems, corresponds to the set of inference rules that is known as the rules for identity types in intensional Martin-Lof type theories (first published in \cite{MLTT73})\footnote{There is also a simpler set of rules corresponding to the identity types in the extensional Martin-Lof type theory (first published in \cite{MLTT79}). Cartmell, in his notion of a strong M-L structure \cite[p.3.36]{Cartmell0}, considers the set of rules for the extensional theory.}. Since the key ingredient of this structure is known in type theory as the J-eliminator we call it the J-structure.

We do not use the ``sequent'' notation that is so widespread in the literature on type theory for general C-systems restricting its use only to examples where we assume the C-system to be a syntactic one.

The reason for this restriction is that translating the sequent-like notations into the algebraic notation of C-systems or categories with families requires considerable mastery of various conventions connected to the use of dependently typed systems. An example of such a translation is the description of an object $IdxT(T)$ corresponding to the sequent-like expression $(\Gamma, x:T, y:T, e:IdT\, T\, x\, y\,;)$ in Construction \ref{2015.03.27.constr1}. 

Some of the difficulties that arise here can already be seen on the translation of the sequent-like expression $(\Gamma, x:T, y:T;)$. Here the same letter $T$ is used to refer to objects of two different types - the first $T$ refers to an element on $Ob_1(\Gamma)$ and the second $T$ refers to an element in $Ob_1(T)$. It is ``understood'' that the second $T$ is the image of the first $T$ under the map $p_T^*:Ob_1(\Gamma)\sr Ob_1(T)$ but this understanding is a part of a tradition and  is not reflected in any mathematical statement that one can refer to. 

For the syntactic C-systems we are allowed to use the sequent notation for the following reason. First, since $CC$ in this case is a sub-quotient of $CC(R,LM)$ our notation only needs to provide a reference to elements of sets associated with $CC(R,LM)$ itself.  There, the first $T$ refers to an element of $LM(\{1,\dots,l\})$ where $l$ is the length of $\Gamma$ and $LM(X)$ is the set of type expressions in the raw syntax with free variables from a set $X$ and the second $T$ refers to an element of $LM(\{1,\dots,l+1\})$ that is the image of the first $T$ under the map
$$LM(\{1,\dots,l\})\sr LM(\{1,\dots,l+1\})$$
defined by the inclusion $\{1,\dots,l\}\subset \{1,\dots,l+1\}$. In this case the map does not depend on $T$. We should distinguish between $IdT$ as a structure on the C-system and the corresponding syntactic construction (because they have different types). If we denote the syntactic ``identity types'' by $IdT^s\, T\, t_1\, t_2$ then for the sequence 
$$\Gamma, x:T, y:T, e:IdT^s\,T\,x\,y;$$
to define an element of $Ob(CC(LM,R))$, the expression $IdT^s\,T\,x\,y$ must refer to an element of $LM(\{1,\dots,l+2\})$ and its form shows that we assume that there is an operation
$$IdT^s:LM\times R\times R\sr R$$
(a natural transformation of functors that is a linear morphism of left $R$-modules) and $IdT^s\,T\,x\,y$ is the ``named variables'' notation for $IdT^s_{{1,\dots,l+2}}(T,l+1,l+2)$. 

We do not continue this explanation of how to construct J-structures on syntactic C-systems. This will be done in a separate paper. Let us remark however that constructing J-structures on syntactic C-systems is relatively easy and that the difficult questions about J-structures on such C-systems are the ones related to the initially properties of the resulting objects. 

While constructing J-structures on the syntactic C-systems relatively straightforward, constructing non-degenerate\footnote{See Remark \ref{2015.05.12.rem1}.} J-structures on semantic C-systems or categories with families proved to be very difficult.

The first instance of such a construction, due to Martin Hofmann and Thomas Streicher, appeared in \cite{Hofmann1}. It was done in the language of categories with families and the underlying category there was the category of groupoids. 

The construction of Hofmann and Streicher was substantially extended and generalized in the Ph.D. thesis of Michael Warren \cite{WarrenThesisProsp},\cite{WarrenThesis} and his subsequent papers such as \cite{Warreninfty}. 

Further important advances were achieved in the work of Richard Garner and Benno van den Berg \cite{BergandGarner}. 

\comment{However the original expectation that it should be possible to construct C-systems or categories with families with J-structures corresponding to all Quillen closed model categories with sufficiently good properties have not been realized. In particular none of the previous methods provided a construction of a C-system whose underlying category is the category of simplicial sets and whose J-structure corresponds, in an appropriately defined sense,  to the standard closed model structure on this category.

This goal is still not fully realized in this paper since to achieve it one has to construct a Kan fibration with certain properties and discussing such a construction is outside of the scope of this paper. }

Two main results of the first part of this paper provide a new approach to the construction of J-structures on semantic C-systems, an approach that can be used to construct the J-structure on the C-system of the univalent model. 

Construction \ref{2015.05.22.constr1} provides a simple way of extending a J1-structure on a universe $p$ in a category $\cal C$ to a full J-structure. 

Construction \ref{2015.04.04.constr2} provides a method of constructing a J-structure on the C-system $CC({\mathcal C},p)$ from a J-structure on $p$.

Combined together they provide a method of constructing a J-structure on $CC({\cal C},p)$ from a J1-structure on $p$.

We also discuss two sets of conditions on a pair of classes of morphisms $TC$ and $FB$ in a locally cartesian closed category that can be used in combination with Construction \ref{2015.05.22.constr1} to construct J-structures. These conditions often hold for the classes of trivial cofibrations and fibrations in model categories (or categories with weak factorization systems) providing a way of constructing C-systems with J-structures starting from such categories. 

In this paper we continue to use the diagrammatic order of writing composition of morphisms, i.e., for $f:X\sr Y$ and $g:Y\sr Z$ the composition of $f$ and $g$ is denoted by $f\circ g$.

In this paper, as in the preceding papers \cite{Cfromauniverse} and \cite{fromunivwithPi}, we often have to consider structures on categories that are not invariant under equivalences and their interaction with structures that are invariant under the equivalences. 

The methods used in this paper are fully constructive and the paper is written in ``formalization ready'' style with all the proofs provided in enough detail to ensure that there are no hidden difficulties for the formalization of all of the results presented here.

Except for the section that discusses the use of classes $TC$ and $FB$, the methods we use are also completely elementary in the sense that they rely only on the essentially algebraic language of categories with various structures. 

The key Definition \ref{2015.03.27.def6} and its relation to the J-structures on categories $CC({\cal C},p)$ first appeared in \cite{CMUtalk}. 

The author would like to thank the Department of Computer Science and Engineering of the University of Gothenburg and Chalmers University of Technology for its hospitality during the work on this paper.

\section{J-structures on C-systems and on universe categories}

\subsection{The J-structure on a C-system}

To define the J-structure on a C-system we will actually define three structures J0-structure, J1-structure over a J0-structure and and J2-structure over a J1-structure with the J-structure being the same as a triple $(IdT,refl,J)$ where $Idt$ is a J0-structure, $refl$ is a J1-structure over $IdT$ and $J$ is a J2-structure over $refl$. For the notations used below see \cite{fromunivwithPi}.
\begin{definition}
\llabel{2015.03.27.def1}
A J0-structure on a C-system $CC$ is a family of functions 
$$IdT_{\Gamma}:\{o_1,o_2\in \wt{Ob}_1(\Gamma)\,|\,\partial(o_1)=\partial(o_2)\}\sr Ob_1(\Gamma)$$
given for all $\Gamma\in Ob$ that is natural in $\Gamma$ i.e. such that for any $f:\Gamma'\sr \Gamma$ and $o,o'\in \wt{Ob}_1(\Gamma)$ with $\partial(o)=\partial(o')$, one has $f^*(IdT_{\Gamma}(o,o'))=IdT_{\Gamma'}(f^*(o),f^*(o'))$.
\end{definition}
\begin{definition}
\llabel{2015.03.27.def2}
Let $IdT$ be a J0-structure on $CC$. A J1-structure over $IdT$ is a family of functions
$$refl:\wt{Ob}_1(\Gamma)\sr \wt{Ob}_1(\Gamma)$$
given for all $\Gamma\in Ob$ such that:
\begin{enumerate}
\item $refl$ is natural in $\Gamma$,
\item for any $\Gamma$ and $o\in \wt{Ob}_1(\Gamma)$ one has 
\begin{eq}
\llabel{2015.03.27.eq8}
\partial(refl(o))=IdT(o,o)
\end{eq}
\end{enumerate}
\end{definition}
To define the notion of a J2-structure over a given J1-structure we will need to describe two constructions first.
\begin{problem}
\llabel{2015.03.27.prob1}
Given a J0-structure $IdT$ to construct a family of functions
$$IdxT:Ob_1(\Gamma)\sr Ob_3(\Gamma)$$
such that for $f:\Gamma'\sr \Gamma$ and $T\in Ob_1(\Gamma)$ one has $f^*(IdxT(T))=IdxT(f^*(T))$.
\end{problem}
\begin{construction}
\llabel{2015.03.27.constr1}\rm
Recall that for $T\in Ob_1(\Gamma)$ we let $\delta(T)$ denote the morphism $T\sr p_T^*(T)$ that can be described equivalently as $s_{Id_T}$ or as the only morphism such that $\delta(T)\circ p_{p_T^*(T)}=Id_T$ and $\delta_T\circ q(p_T,T)=Id_T$. Because of the first equation we have $\delta(T)\in \wt{Ob}(p_T^*(T))$. 

We define:
\begin{eq}
\llabel{2015.04.06.eq1}
IdxT(T):=IdT_{p_T^*(T)}((p_{p_T^*(T)}^*(\delta(T))), \delta(p_T^*(T)))
\end{eq}
We have
$$p_{p_T^*(T)}^*(\delta_T)\in \wt{Ob}(p_{p_T^*(T)}^*(p_T^*(T)))$$
and 
$$\delta(p_T^*(T))\in \wt{Ob}(p_{p_T^*(T)}^*(p_T^*(T)))$$
and since $ft(p_{p_T^*(T)}^*(p_T^*(T)))=p_T^*(T)$ the $IdxT(T)$ is well defined. The fact that $IdxT(T)\in Ob_3(\Gamma)$ follows now from the fact that $ft^2(p_T^*(T))=ft(T)=\Gamma$. The objects and some of the morphisms involved in this construction can be seen on the diagram: 
$$
\begin{CD}
p_{p_T^*(T)}^*(p_T^*(T)) @>>> p_T^*(T) @>>> T\\
@VVV @Vp_{p_T^*(T)} VV @VV p_TV\\
p_T^*(T) @>p_{p_T^*(T)}>> T @>p_T>> \Gamma
\end{CD}
$$
The proof that $IdxT$ is natural in $f:\Gamma'\sr \Gamma$ is straightforward.
\end{construction}
\begin{problem}
\llabel{2015.03.27.prob2}
Given a J0-structure $IdT$ and a J1-structure $refl$ over it to construct for all $\Gamma\in Ob$ and $T\in Ob_1(\Gamma)$ a morphism 
$$rf_T:T\sr IdxT(T)$$
such that for any $f:\Gamma'\sr \Gamma$ one has $f^*(rf_T)=rf_{f^*(T)}$. 
\end{problem}
\begin{construction}
\llabel{2015.03.27.constr2}\rm
We have:
$$\delta(T)^*(IdxT(T))=\delta(T)^*(IdT_{p_T^*(T)}((p_{p_T^*(T)}^*(\delta_T))),\delta(p_T^*(T)))$$
$$=IdT_T(\delta(T)^*(p_{p_T^*(T)}^*(\delta_T)),\delta(T)^*(\delta(p_T^*(T))))=$$
$$=IdT_T(\delta(T),\delta(T))$$
where the last equality follows from the fact that $f^*(\delta(T))=s_f$ and for any $s\in\wt{Ob}$, $s_s=s$.  This shows that we have a canonical square of the form
\begin{eq}
\llabel{2015.03.31.eq3}
\begin{CD}
IdT(\delta(T),\delta(T)) @>q(\delta(T),IdxT(T))>> IdxT(T)\\
@VVV @VVV\\
T @>\delta(T)>> p_T^*(T)
\end{CD}
\end{eq}
and $refl(\delta(T))$ is a morphism $T\sr IdT(\delta(T),\delta(T))$. We define:
\begin{eq}
\llabel{2015.04.02.eq1}
rf_T:=refl(\delta(T))\circ q(\delta(T),IdxT(T))
\end{eq}
The proof that for any $f:\Gamma'\sr \Gamma$ one has $f^*(rf_T)=rf_{f^*(T)}$ is straightforward. 
\end{construction}
\begin{definition}
\llabel{2015.03.27.def3}
Let $IdT$ and $refl$ be a J0-structure and a J1-structure over it. A J2-structure over $(IdT,refl)$ is data of the form: for all $\Gamma\in Ob$, for all $T\in Ob_1(\Gamma)$, for all $P\in Ob_1(IdxT(T))$, for all $s0\in \wt{Ob}(rf_T^*(P))$, an element $J(\Gamma,T,P,s0)$ of $\wt{Ob}(P)$ such that:
\begin{enumerate}
\item $J$ is natural in $\Gamma$, i.e., for any $f:\Gamma'\sr \Gamma$ and $T,P,s0$ as above one has
$$f^*(J(\Gamma,T,P,s0))=J(\Gamma',f^*(T),f^*(P),f^*(s0))$$
where the right hand side of the equation is well-defined because of the naturality in $f$ of $IdxT$ and $rf$. 
\item $J$ satisfies the $\iota$-rule. For $\Gamma, T, P, s0$ as above one has
$$rf_T^*(J(\Gamma,T,P,s0))=s0$$
\end{enumerate}
\end{definition}
\begin{remark}
\llabel{2015.05.12.rem1}\rm
A J0-structure is called degenerate or extensional if for all $T\in Ob_{\ge 1}(CC)$ and $o,o'\in \wt{Ob}(T)$ one has\footnote{The following is the classical way of saying that there is an equivalence between the types $\wt{Ob}(IdxT(o,o'))$ and $(o=o')$.}
$$\wt{Ob}(IdxT(o,o'))=\left\{
\begin{array}{ll}
\emptyset&\mbox{\rm if $o\ne o'$}\\
pt&\mbox{\rm if $o=o'$}
\end{array}
\right.
$$
One can easily see that any two extensional J0-structures are equal and that any extensional J0-structure has a unique extension to a full J-structure that is also called extensional. 

We will not consider these extended versions of $J$ in the present version of the paper. 
\end{remark}

\begin{remark}\rm
\llabel{2015.05.24.rem1}
When one studies J-structures on C-systems that have no $(\Pi,\lambda)$-structures it is important, as emphasized for example in \cite{vandenBergGarner2011}, to consider a more complex structure than the one that we consider here. This more complex structure can be seen as a family of structures $eJ_n$ where $eJ_0=J2$ and where $eJ_n$ over $(IdT,refl)$ is a collection of data of the form: for all $\Gamma\in Ob$, for all $T\in Ob_1(\Gamma)$, for all $\Delta\in Ob_{n}(IdxT(T))$, for all $P\in Ob_1(\Delta)$, for all $s0\in \wt{Ob}(rf_T^*(P))$, an element $eJ_n(\Gamma,T,\Delta,s0)$ in $\wt{Ob}(P)$ such that $eJ_n$ satisfies the obvious analog of $\iota$-rule and such that it is natural in $\Gamma$. See also Remark \ref{2015.05.24.rem2}. 
\end{remark}

\subsection{The J-structure on a universe in a category}
Let $\mathcal C$ be a (pre-)category and $p:\wt{U}\sr U$ a morphism in $\mathcal C$. Recall that a universe structure on $p$ is a choice of pull-back squares of the form
$$
\begin{CD}
(X;F) @>Q(F)>> \wt{U}\\
@Vp_{X,F} VV @VV p V\\
X @>F>> U
\end{CD}
$$
for all $X$ and all morphisms $F:X\sr U$. A universe in $\mathcal C$ is a morphism with a universe structure on it and a universe category is a category with a universe and a choice of a final object $pt$.  

For $f:W\sr X$ and $g:W\sr \wt{U}$ we will denote by $f*g$ the unique morphism such that 
$$(f*g)\circ p_{X,F}=f$$
$$(f*g)\circ Q(F)=g$$
When we need to distinguish canonical squares of different universes we may write $(X;F)_{p}$ and $f*_p g$.
\begin{remark}\rm
\llabel{2015.03.29.rm1}
Note that we made no assumption about $Q(Id_U)$ being equal to $Id_{\wt{U}}$. In fact, since we want the results of this paper to be constructive, we are not allowed to make such an assumption since the question of whether or not a given morphism is an identity morphism need not be decidable and therefore we can not ``normalize'' our constructions by doing a ``case'' on whether a morphism is an identity morphism or not. The importance of this observation (in the context of whether a simplex is degenerate or not) was emphasized by \cite{BCH}.
\end{remark}
For $X'\stackrel{f}{\sr}X\stackrel{F}{\sr}U$ we let $Q(f,F)$ denote the morphism 
$$(p_{X',f\circ F}\circ f)*Q(f\circ F):(X';f\circ F)\sr (X;F)$$
As is shown in \cite{fromunivwithPi}, the square
\begin{eq}
\llabel{2015.04.06.l0.sq}
\begin{CD}
(X';f\circ F) @>Q(f,F)>> (X;F)\\
@Vp_{X',f\circ F} VV @VVp_{X,F}V\\
X' @>f>> X
\end{CD}
\end{eq}
is a pull-back square.

Following \cite{fromunivwithPi} we define for any universe $p:\wt{U}\sr U$ and any $V\in {\mathcal C}$ a functor
$$D_p(-,V):X\mapsto \amalg_{F:X\sr U}Hom((X;F), V)$$
whose action on morphisms is given by
$$D_p(f,V):(F,a)\mapsto (f\circ F, Q(f,F)\circ a)$$
When $\mathcal C$ is a locally cartesian closed category any morphism $p:\wt{U}\sr U$ defines a functor
$$I_p:V\mapsto \uu{Hom}((\wt{U},p),(U\times V,pr_1))$$
and we have constructed in \cite[Construction 3.9]{fromunivwithPi} a family of bijections
$$\eta^!_{p,X,V}:Hom(X,I_p(V))\sr D_p(X,V)$$
that are natural in $X$ and $V$. We let $\eta$ denote the inverse bijections 
$$\eta_{p,X,V}:D_p(X,V)\sr Hom(X,I_p(V))$$
Using the functorial structure on the mapping $V\mapsto (U\times V,pr_1)$ together with the naturality of internal Hom-objects in the second argument we get a functoriality structure on $I_p$
$$(f:V\sr V')\mapsto (I_p(f):I_p(V)\sr I_p(V'))$$
Similarly, using the functoriality of $\uu{Hom}$ in the second argument (see e.g. the appendix in \cite{fromunivwithPi}) we obtain, for any $p:\wt{U}\sr U$, $p':\wt{U}'\sr U$ and $h:\wt{U}'\sr \wt{U}$ over $U$ and $V$ a morphism
$$I^h(V):I_p(V)\sr I_{p'}(V)$$
\begin{lemma}
\llabel{2015.04.10.l2}
In the notations introduced above let $f:V\sr V'$ be a morphism, then the square
$$
\begin{CD}
I_{p'}(V) @>I_{p'}(f)>> I_{p'}(V')\\
@V{I^h(V)}VV @VV{I^h(V')}V\\
I_p(V) @>I_p(f)>> I_p(V')
\end{CD}
$$
\end{lemma}
\begin{myproof}
This is a particular case of the commutative square of \cite[Lemma 8.5]{fromunivwithPi}.
\end{myproof}

\begin{lemma}
\llabel{2015.04.02.l4}
Let $p:\wt{U}\sr U$ and $p':\wt{U}'\sr U$ be two morphisms with universe structures and $f:\wt{U}'\sr \wt{U}$ be a morphism over $U$. For $V\in {\mathcal C}$ let $I^f(V)$ be the corresponding morphism $I_{p'}(V)\sr I_p(V)$. Then for any $X$ the square
$$
\begin{CD}
D_p(X,V) @>\eta_{p,X,V}>> Hom(X, I_p(V))\\
@V{D^f(X,V)}VV @VV{-\circ I^f(V)} V\\
D_{p'}(X,V) @>\eta_{p',X.V}>> Hom (X, I_{p'}(V))
\end{CD}
$$
where the left hand side arrow is of the form 
$$D^f(X,V):(F, F')\mapsto (F,F^*(f)\circ F')$$
commutes.
\end{lemma}
\begin{myproof}
Since $\eta$ is defined as an inverse to $\eta^!$ it is sufficient to show that for any $g\in Hom(X,I_p(V))$ one has $\eta^{',!}(g\circ I^f(V))=D^f(X,V)(\eta^!(g))$. Let
$$pr=prI_p(V):I_p(V)\sr U$$
$$pr'=prl_{p'}(V):I_{p'}(V)\sr U$$
be the canonical projections. Let
$$st=st_p(V):(I_p(V);pr)\sr V$$
$$st'=st_{p'}(V):(I_{p'}(V);pr')'\sr V$$
be the morphisms introduced in \cite{fromunivwithPi}. By \cite[Problem 3.8]{fromunivwithPi} we have
$$\eta^{',!}(g\circ I^f(V))=(g\circ I^f(V)\circ pr', Q'(g\circ I^f(V), pr')\circ st')$$
and
$$D^f(X,V)(\eta^!(g))=D^f(X,V)(g\circ pr, Q(g,pr)\circ st)=(g\circ pr, (g\circ pr)^*(f)\circ Q(g,pr)\circ st)$$
Therefore it is sufficient to show that
$$I^f(V)\circ pr'=pr$$
and
$$Q'(g\circ I^f(V), pr')\circ st'=(g\circ pr)^*(f)\circ Q(g,pr)\circ st$$
The first equality asserts that $I^f(V)$ is a morphism over $U$ which follows from its construction. 

By Lemma \ref{2015.04.20.l1} we have 
$$(g\circ pr, (g\circ pr)^*(f)\circ Q(g,pr)=Q'(g,pr)\circ pr^*(f)$$
Next we have
$$Q'(g\circ I^f(V), pr')=Q'(g,I^f(V)\circ pr')\circ Q'(I^f(V),pr')=Q'(g,pr)\circ Q'(I^f(V),pr')$$
by \cite[Lemma 3.2]{fromunivwithPi}. It remains to check that
$$Q'(I^f(V),pr')\circ st'=pr^*(f)\circ st$$
This requires opening up the definitions of $st$ and $st'$ which gives us
$$Q'(I^f(V),pr')\circ \iota'\circ ev'\circ pr_2=pr^*(f)\circ \iota\circ ev\circ pr_2$$
We will obtain this equality as a consequence of commutativity of three squares:
$$
\begin{CD}
(I_p(V);pr)' @>Q'(I^f(V),pr)>> (I_{p'}(V);pr')'\\
@V\iota' VV @V\iota' VV\\
(I_p(V),pr)\times_U(\wt{U}',p') @>I^f(V)\times Id>> (I_{p'}(V),pr')\times_U(\wt{U}',p')
\end{CD}
$$
$$
\begin{CD}
(I_p(V);pr)'  @>pr^*(f)>> (I_p(V);pr)\\
@V\iota' VV @VV\iota V\\
(I_p(V),pr)\times_U(\wt{U}',p') @>Id\times f >> (I_p(V),pr)\times_U(\wt{U},p)
\end{CD}
$$
and
$$
\begin{CD}
(I_p(V),pr)\times_U(\wt{U}',p') @>I^f(V)\times Id >> (I_{p'}(V),pr')\times_U(\wt{U}',p')\\
@VId\times f VV @VV {ev'}V\\
(I_p(V),pr)\times_U(\wt{U},p) @>ev>> U\times V
\end{CD}
$$
The first two squares are particular cases of \cite[Lemma 8.1]{fromunivwithPi}. To obtain the first one one has to set $Z=U$, $b=Id_{\wt{U}'}$, and $a=I^f(V)$. To obtain the second one one has to set $Z=U$, $b=f$ and $a=Id_{I_p(V)}$. The last square is a particular case of \cite[Lemma 8.6]{fromunivwithPi}.
%
\end{myproof}
\begin{definition}
\llabel{2015.03.27.def4}
A J0-structure on a universe $p$ in a category $\mathcal C$ is a morphism $Eq:(\wt{U};p)\sr U$.
\end{definition}

Let $Eq$ be a J0-structure on $p$. Consider the object 
$$E\wt{U}:=(\wt{U};p,Eq)$$
as an object over $U$ relative to the composition of projections 
$$E\wt{U}\stackrel{p_{(\wt{U};p),Eq}}{\lr}(\wt{U};p)\stackrel{p_{\wt{U},p}}{\lr} \wt{U}\stackrel{p}{\sr} U$$
that we will denote by $pE\wt{U}$. 
\begin{problem}
\llabel{2015.05.08.prob1}
To construct a universe structure on $pE\wt{U}$.
\end{problem}
\begin{construction}\rm
\llabel{2015.05.08.constr1}
We have three diagrams with pull-back squares of the form:
$$
\begin{CD}
(X;F,Q(F)\circ p, Q(Q(F),p)\circ Eq) @>Q(Q(Q(F),p),Eq)>> (\wt{U}; p, Eq)\\
@VVV @VVV\\
(X;F,Q(F)\circ p) @>Q(Q(F),p)>> (\wt{U};p)
\end{CD}
$$$$
\begin{CD}
(X;F,Q(F)\circ p) @>Q(Q(F),p)>> (\wt{U},p)\\
@VVV @VVp_{\wt{U},p}V\\
(X,F) @>Q(F)>> \wt{U}
\end{CD}
$$$$
\begin{CD}
(X;F) @>Q(F)>> \wt{U}\\
@Vp_{X,F} VV @VVV\\
X @>F>> U
\end{CD}
$$
and we define the canonical square for $F$ relative to $pE\wt{U}$ to be the square obtained by concatenating these three squares vertically. 
\end{construction}
Let us denote the components of the canonical squares for $pE\wt{U}$ as follows:
$$
\begin{CD}
(X;F)_{E} @>Q(F)_{E}>> E\wt{U}\\
@Vp_{X,F}^{E} VV @VVpE\wt{U} V\\
X @>F>> U
\end{CD}
$$
Explicitly we have
$$(X;F)_{E}=(X;F,Q(F)\circ p, Q(Q(F),p)\circ Eq)$$
$$Q(F)_{E}=Q(Q(Q(F),p),Eq)$$
$$p_{X,F}^E=p_{(X;F,Q(F)\circ p),Q(Q(F),p)\circ Eq}\circ p_{(X;F),Q(F)\circ p}\circ p_{X,F}$$
We will also write $Q(f,F)_{E}$ for the canonical morphisms from $(X;f\circ F)_{E}$ to $(X;F)_{E}$. 
\begin{definition}
\llabel{2015.03.27.def5}
Let $p$ be a universe in $\mathcal C$ and $Eq$ be a J0-structure on $p$. A J1-structure on $p$ over $Eq$ is a morphism $\Omega:\wt{U}\sr \wt{U}$ such that the square
\begin{eq}\llabel{2015.03.27.sq1}
\begin{CD}
\wt{U} @>\Omega>> \wt{U}\\
@V\Delta VV @VVp V\\
(\wt{U};p) @>Eq >> U
\end{CD}
\end{eq}
where $\Delta=(Id_{\wt{U}})*(Id_{\wt{U}})$ is the diagonal of $\wt{U}$, commutes.
\end{definition}

The square (\ref{2015.03.27.sq1}) defines a morphism $\wt{U}\sr E\wt{U}$ that we will denote by $\omega$.

To define a J2-structure on a universe we will need to assume that $\mathcal C$ is a locally cartesian closed category.  Recall that locally cartesian closed category is a category with the choice of fiber squares based on all pairs of morphisms with a common codomain as well as the choice of relative internal Hom-objects and co-evaluation morphisms for all such pairs. For our notations related to the locally cartesian closed categories as well as for some other notations used below see \cite{fromunivwithPi}. 

When a universe is considered in a locally cartesian closed category we make no assumption about the compatibility of choices of the pull-back squares of the universe structure on $p$ and pull-back squares of the locally cartesian closed structure. 

Consider the functors $I_{p}$ and $I_{pE\wt{U}}$. We have the following commutative square:
\begin{eq}\llabel{2010.sq1}
\begin{CD}
I_{pE\wt{U}}(\wt{U}) @>I^{\omega}(\wt{U})>> I_p(\wt{U})\\
@VI_{pE\wt{U}}(p)VV @VVI_p(p)V\\
I_{pE\wt{U}}(U) @>I^{\omega}(U)>> I_p(U)
\end{CD}
\end{eq}
and therefore a morphism
$$I_{pE\wt{U}}(\wt{U}) \stackrel{coJ}{\lr}
(I_{pE\wt{U}}(U), I^{\omega}(U)) \times_{I_p(U)} (I_p(\wt{U}), I_p(p))
$$
\begin{definition}
\llabel{2015.03.27.def6}
A J2-structure on $p$ relative to a J0-structure $Eq$ and J1-structure $\Omega$, is a morphism 
$$
Jp:( I_{pE\wt{U}}(U), I^{\omega}(U))\times_{I_p(U)} (I_p(\wt{U}), I_p(p))\sr I_{pE\wt{U}}(\wt{U}) $$
such that $Jp\circ coJ = Id$.
\end{definition}
Note that we have:
\begin{eq}
\llabel{2015.04.04.eq1}
J\circ I^{\omega}(\wt{U})=J\circ coJ\circ prI_p(\wt{U})=prI_{p}(\wt{U})
\end{eq}
\begin{eq}
\llabel{2015.04.04.eq2}
J\circ I_{pE\wt{U}}(p)=J\circ coJ\circ prI_{pE\wt{U}}(U)=prI_{pE\wt{U}}(U).
\end{eq}
where $prI_p(V)$ is the canonical morphism $I_p(V)\sr U$.

A J-structure on $p$ is a triple $(Eq,\Omega,Jp)$ where $Eq$ is a J0-structure, $\Omega$ is a J1-structure relative to $Eq$ and $Jp$ a J2-structure relative to $Eq$ and $\Omega$.

For a J1-structure $(Eq,\Omega)$ on a universe in a category with a locally cartesian closed structure let $Fp_{Eq,\Omega}$ denote the fiber product
$$(I_{pE\wt{U}}(U), I^{\omega}(U)) \times_{I_p(U)} (I_p(\wt{U}), I_p(p))$$
and let $pFp_{Eq,\Omega}=I^{\omega}(U)\dd I_p(p)$ be the projection $Fp_{Eq,\Omega}\sr U$. Let further $pr_1$ be the projection from $Fp$ to $I_{pE\wt{U}}(U)$ and $pr_2$ the projection from $Fp$ to $I_p(\wt{U})$.

Our solution to the following problem is the key to the construction of J-structures over a given J1-structure in categories with weak factorization systems in particular in Quillen model categories. 
\begin{problem}
\llabel{2015.05.12.l1}
Let $\mathcal C$ be a category with a locally cartesian closed structure and $Eq,\Omega$ be a J1-structure on $({\mathcal C},p)$. To construct a bijection 
between the set of J-structures on $p$ over $(Eq,\Omega)$ and the set of morphisms $(Fp,pFp)\times_U(E\wt{U},pE\wt{U})\sr \wt{U}$ that split the following square into two commutative triangles:
\begin{eq}\llabel{2015.05.22.sq1}
\begin{CD}
(Fp,pFp)\times_U(\wt{U},p) @>adj(pr_2)\circ pr_2>> \wt{U}\\
@VId_{Fp}\times \omega VV @VV p V\\
(Fp,pFp)\times_U(E\wt{U},pE\wt{U}) @>adj(pr_1)\circ pr_2>> U
\end{CD}
\end{eq}
\end{problem}
\begin{construction}\rm
\llabel{2015.05.22.constr1}
Observe first that there is a bijection between the set of morphisms 
$$(Fp,pFp)\times_U(E\wt{U},pE\wt{U})\sr \wt{U}$$
that split the square (\ref{2015.05.22.sq1}) into two commutative triangles and the set of morphisms 
$$(Fp,pFp)\times_U(E\wt{U},pE\wt{U})\sr U\times\wt{U}$$
that split into two commutative triangles the square:
$$
\begin{CD}
(Fp,pFp)\times_U(\wt{U},p) @>adj(pr_2)>> U\times\wt{U}\\
@VId_{Fp}\times \omega VV @VV Id_U\times p V\\
(Fp,pFp)\times_U(E\wt{U},pE\wt{U}) @>adj(pr_1)>> U\times U
\end{CD}
$$
The rule $f\mapsto adj(f)$ gives us a bijection of the form
$$Hom_U((Fp,pFp),(I_{pE\wt{U}}(\wt{U}),\_))\sr Hom_U((Fp,pFp)\times_U (E\wt{U},pE\wt{U}), (U\times\wt{U}, pr_2))$$
All sections of $coJ$ are automatically morphisms over $U$. Therefore it remains to show that this bijection defines a bijection of the subset of morphisms that are sections of $coJ$ and morphisms that make the two triangles commutative. 

One verifies first that a morphism $f:Fp\sr I_{pE\wt{U}}(\wt{U})$ is a section of $coJ$ if and only if $f\circ I_{pE\wt{U}}(p)=pr_1$ and $f\circ I^{\omega}(\wt{U})=pr_2$. This is straightforward. 

Next we have
$$I_{pE\wt{U}}(p)=\uu{Hom}_U((E\wt{U},pE\wt{U}),Id_U\times p)$$
$$I^{\omega}(\wt{U})=\uu{Hom}_U(\omega,(U\times\wt{U},pr_2))$$
Therefore by \cite[Lemma 8.7]{fromunivwithPi} one has
$$adj(f\circ I_{pE\wt{U}}(p))=adj(f)\circ (Id_U\times p)$$
$$adj(f\circ I^{\omega}(\wt{U}))=(Id_{Fp}\times_{U}\omega)\circ adj(f)$$
and we conclude that $f$ is a section of $coJ$ iff 
$$adj(f)\circ (Id_U\times p)=adj(pr_1)$$
$$(Id_{Fp}\times_{U}\omega)\circ adj(f)=adj(pr_2)$$
This completes the construction.
\end{construction}
\begin{remark}\rm
\llabel{2015.05.24.rem2}
It is likely to be relatively easy to generalize the constructions of this  paper to the extended J-structures $eJ_n$ (see Remark \ref{2015.05.24.rem1}). The key to such generalization is \cite[Remark 3.13]{fromunivwithPi}. The structures $eJp_n$ can be defined in the same way as $Jp$ but with the square (\ref{2010.sq1}) replaced by the square
\begin{eq}
\llabel{2015.05.24.sq1}
\begin{CD}
I_{pE\wt{U}}(I_p^n(\wt{U})) @>I^{\omega}(I_p^n(\wt{U}))>> I_p(I_p^n(\wt{U}))\\
@VI_{pE\wt{U}}(I_p^n(p)) VV @VV I_p(I_p^n(p))V\\
I_{pE\wt{U}}(I_p^n(U)) @>I^{\omega}(I_p^n(U))>> I_p(I_p^n(U))
\end{CD}
\end{eq}
\end{remark}

\subsection{J-structures on universes in categories with two classes of morphisms}
This is the only part of the paper where we depart from constructions that are conservatively algebraic over the theory of categories, i.e., from constructions that can be expressed in terms of adding new essentially algebraic operations to the theory of categories without adding new sorts to this theory. 

Considering classes of morphisms in categories can be expressed in the essentially algebraic way but this requires adding new sorts to the theory. 

This is also the only context where we use the concept ``there exists'' in this paper. In all the previous cases the objects that we considered were given (specified). To make the lemmas proved below into constructions and to avoid the use of ``there exists'' one would have to define the collection $FB$ as a collection of pairs of a morphism $p$ together with, for all $i\in TC$, $f_W$ and $f_Z$ such that $f_Z\circ p=i\circ f_W$, a morphism $g$ such that $i\circ g=f_W$ and $g\circ p=f_W$. 

Recall that a collection of morphisms $R$ is said to have the right lifting property for the collection of morphisms $L$ if for any commutative square of the form
$$
\begin{CD}
Z @>f_Z>> E\\
@Vi VV @VVp V\\
W @>f_W>> B
\end{CD}
$$
such that $i\in L$ and $p\in R$ there exists a morphism $g:W\sr E$ that makes the two triangles into which it splits the square to commute i.e. a morphism $g$ such that $i\circ g=f_Z$ and $g\circ p=f_W$. 

We are going to consider two sets of conditions (Conditions \ref{2015.05.22.cond2} and \ref{2015.05.22.cond1}) on a pair of classes of morphisms $FB$ and $TC$ in a category with fiber products and then show in Theorems \ref{2015.05.22.th1} and \ref{2015.05.16.th1} how pairs satisfying conditions of each of these two sets can be used to construct J-structures on elements of $FB$. 

Our first set of conditions is as follows:
\begin{cond}\llabel{2015.05.22.cond2}
\begin{enumerate}
\item A morphism is in $FB$ if and only if it has the right lifting property for $TC$,
\item consider morphisms $f: B'\sr B$, $p_1:E_1\sr B$, $p_2:E_2\sr B$ and $i:E_1\sr E_2$ such that $p_1,p_2\in FB$ and $i\in TC$. Then the morphism
$$Id_{B'}\times i: (B',f)\times_B(E_1,p_1)\sr (B',f)\times_B(E_2,p_2)$$
is in $TC$.
\end{enumerate}
\end{cond}
\begin{theorem}
\llabel{2015.05.22.th1}
Let $FB$ and $TC$ be two classes of morphisms in a locally cartesian closed category $\cal C$ that satisfy Conditions \ref{2015.05.22.cond2}. Let $p$ be a universe in $\cal C$ and $(Eq,\Omega)$ a J1-structure on $p$ such that:
\begin{enumerate}
\item $p$ is in $FB$,
\item $\omega$ is in $TC$.
\end{enumerate}
Then there exists an extension of $(Eq,\Omega)$ to a full J-structure on $p$.
\end{theorem}
\begin{myproof}
Let us apply Construction \ref{2015.05.22.constr1} to $(Eq,\Omega)$. To construct the required morphism it is sufficient to establish that $Id_{Fp}\times\omega$ is in $TC$. It follows from the first of our conditions that $FB$ is closed under pull-backs and compositions. Therefore, $pE\wt{U}$ is in $FB$. It remains to apply the second of our conditions.
\end{myproof}

Our second set of conditions is more involved. Conditions of this set can be satisfied in the situations arising when one attempts to localize Quillen model structures and when the resulting sets of morphisms do not for a model structure. The difference is mainly concerned with the fact that the good behavior is required only for fibrations over fibrant objects. One particular example of such situation is considered in \cite[Section 3.3]{SRF}.

\begin{cond}\llabel{2015.05.22.cond1}
\begin{enumerate}
\item $Id_{pt}$ is in $FB$,
\item let $B$ be such that the morphism $B\sr pt$ is in $FB$ then a morphism $p:E\sr B$ is in $FB$ if and only if it has the right lifting property for $TC$,
\item if $p:E\sr B$ and $B\sr pt$ are in $FB$, $i:Z\sr W$ is in $TC$ and $f:W\sr B$ iis an arbitrary morphism then
$$(i\times_U Id_E):(Z,i\circ f)\times_B (E,p)\sr (W,f)\times_B (E,p)$$
is in $TC$.
\end{enumerate}
\end{cond}
We will say that $B$ is fibrant if the morphism $B\sr pt$ is in $FB$. 
\begin{lemma}
\llabel{2015.05.14.l2}
Let $p:E\sr B$ be in $FB$ and $f:B'\sr B$ be a morphism. Assume in addition that $B$ and $B'$ are fibrant, then for any pull-back square of the form
$$
\begin{CD}
E' @>>> E\\
@Vp' VV @VV p V\\
B' @>f>> B
\end{CD}
$$
the morphism $p'$ is in $FB$. 
\end{lemma}
\begin{myproof}
Since $B'$ is fibrant it is sufficient to verify that $p'$ has the right lifting property for $TC$. This can be shown in the standard way to be a consequence of $p$ having the right lifting property for $TC$. That $p$ has this property we know because $p$ is in $FB$ and $B$ is fibrant.
\end{myproof}
\begin{lemma}
\llabel{2015.05.14.l4}
Let $B$ be fibrant and $p_2:E_2\sr E_1$, $p_1:E_1\sr B$ be in $FB$. Then $p_2\circ p_1$ is in $FB$.
\end{lemma}
\begin{myproof}
Let us show first that $E_1$ is fibrant i.e. that $\pi_{E_1}:E_1\sr pt$ is in $FB$. Since $pt$ is fibrant it is sufficient to show that $\pi_{E_1}$ has the right lifting property for $TC$. It is shown in the standard way from the fact that both $p_1$ and $\pi_B:B\sr pt$ have the right lifting property for $TC$ and $\pi_{E_1}=p_1\circ \pi_B$. 

Since $E_1$ is fibrant we know that $p_2$ has the right lifting property for $TC$ and since $B$ is fibrant we know that $p_1$ has the right lifting property for $FB$. From this we conclude in the standard way that $p_2\circ p_1$ have the right lifting property for $TC$ and since $B$ is fibrant this implies that $p_2\circ p_1$ is in $FB$.
\end{myproof}
\begin{lemma}
\llabel{2015.05.14.l1}
Assume that $U,V$ are fibrant and that $p:\wt{U}\sr U$ is in $FB$. Then the morphism $prI_p(V):I_p(V)\sr U$ is in $FB$.
\end{lemma}
\begin{myproof}
Since $U$ is fibrant it is sufficient to check that $pr=prI_p(V)$ has the right lifting property for $TC$. Consider a commutative square of the form
$$
\begin{CD}
Z @>f_Z>> \uu{Hom}_U((\wt{U},p),(U\times V,pr_1))\\
@V i VV @VV pr V\\
W @>f_W >> U
\end{CD}
$$
We need to construct a morphism $f:W\sr \uu{Hom}_U((\wt{U},p),(U\times V,pr_1))$ that would make the two triangles commutative. The commutativity of the lower triangle means that $f$ is a morphism over $U$ which is equivalent to the assumption that $f=adj^{-1}(g)$ for some $g:(W,f_W)\times_U (\wt{U},p)\sr U\times V$ over $U$.

Consider the square
$$
\begin{CD}
(Z, i\circ f_W)\times_U (\wt{U},p) @>adj(f_Z)>> U\times V\\
@V i\times Id_{\wt{U}} VV @VV pr_1 V\\
(W,f_W)\times_U (\wt{U},p) @> f_W\dd p >> U
\end{CD}
$$
By Lemma \ref{2015.05.14.l2} we know that $pr_1$ belongs to $FB$. By our assumptions on $TC$ and $FB$ we know that $i\times Id_{\wt{U}}$ is in $TC$. Therefore there exists a morphism  $g:(W,f_W)\times_U (\wt{U},p) \sr U\times V$ that makes the two triangles commute. 
The commutativity of the lower triangle means that this is a morphism over $U$. Therefore $adj^{-1}(g)$ is defined. Set $f=adj^{-1}(g)$. It remains to check that $i\circ f=f_Z$. This is equivalent to $adj(i\circ f)=adj(f_Z)$. Since $adj(i\circ f)=(i\times Id_{\wt{U}})\circ adj(f)$ by \cite[Lemma 8.7(3)]{fromunivwithPi}, this is equivalent to $(i\times Id_{\wt{U}})\circ g=adj(f_Z)$ which is the commutativity of the upper triangle. 
\end{myproof}
\begin{lemma}
\llabel{2015.05.14.l3}
Assume that $U,V$ are fibrant and that $p:\wt{U}\sr U$ and $r:V'\sr V$ are in $FB$ then $I_p(r):I_p(V')\sr I_p(V)$ is in $FB$.
\end{lemma}
\begin{myproof}
By Lemmas \ref{2015.05.14.l1} and \ref{2015.05.14.l4} we know that $I_p(V)$ is fibrant. Therefore it is sufficient to show that $I_p(r)$ has the right lifting property for $TC$. Consider a commutative square of the form 
\begin{eq}
\llabel{2015.05.14.sq1}
\begin{CD}
Z @>f_Z>> \uu{Hom}_U((\wt{U},p),(U\times V',pr_1))\\
@V i VV @VV \uu{Hom}_U((\wt{U},p), Id_U\times r) V\\
W @>f_W >> \uu{Hom}_U((\wt{U},p), (U\times V, pr_1))
\end{CD}
\end{eq}
The lower right corner is an object over $U$ through the morphism $p\triangle pr_1$. Let $p_W=f_W\circ (p\triangle pr^{U,V}_U)$ and 
$$p_Z=i\circ p_W=f_Z\circ (p\triangle pr^{U,V'}_U)$$
Consider the square
\begin{eq}
\llabel{2015.05.14.sq2}
\begin{CD}
(Z,p_Z)\times_U (\wt{U},p) @>adj(f_Z)>> U\times V'\\
@Vi\times Id_{\wt{U}} VV @VV Id_U\times r V\\
(W,p_W)\times_U (\wt{U},p) @>adj(f_W)>> U\times V
\end{CD}
\end{eq}
This square commutes. Indeed,
$$adj(f_Z)\circ (Id_U\times r )=adj(f_Z\circ \uu{Hom}_U((\wt{U},p), Id_U\times r))=$$
$$adj(i\circ f_W)=(i\times Id_{\wt{U}})\circ adj (f_W)$$
where the first equality is by \cite[Lemma 8.7(1)]{fromunivwithPi} and the third by \cite[Lemma 8.7(3)]{fromunivwithPi}. By Lemmas \ref{2015.05.14.l2} and \ref{2015.05.14.l4} we know that $Id_U\times r$ is in $FB$. By our assumption (3) on $FB$ and $TC$ we know that $i\times Id_{\wt{U}}$ is in $TC$. Therefore, there exists a morphism $g:(W,p_W)\times_U (\wt{U},p) \sr U\times V'$ that splits this square into two commutative triangles. Since the lower triangle commutes, $g$ is a morphism over $U$ and in particular $g=adj(f)$ for some $f:W\sr \uu{Hom}_U((\wt{U},p),(U\times V',pr_1))$. Let us show that $f$ splits the square (\ref{2015.05.14.sq1}) into two commutative triangles i.e. that we have $i\circ f= f_Z$ and $f\circ \uu{Hom}_U((\wt{U},p), Id_U\times r)=f_W$. 

The first equality is equivalent to $adj(i\circ f)=adj(f_Z)$ which is equivalent, by \cite[Lemma 8.7(3)]{fromunivwithPi} to $(i\times Id_{\wt{U}})\circ g=adj(f_Z)$ which is the commutativity of the upper of the two triangles into which $g$ splits (\ref{2015.05.14.sq2}). 

The second equality is equivalent to $adj(f\circ \uu{Hom}_U((\wt{U},p), Id_U\times r))=adj(f_W)$, which is equivalent by \cite[Lemma 8.7(1)]{fromunivwithPi} to 
$g\circ (Id_U\times r)=adj(f_W)$ which is the commutativity of the lower of the two triangles into which $g$ splits (\ref{2015.05.14.sq2}). 

Lemma is proved.
\end{myproof}

We can now prove the second main theorem of this section.
\begin{theorem}
\llabel{2015.05.16.th1}
Let $({\mathcal C},p,pt)$ be a universe category, let $\mathcal C$ be given a locally cartesian closed structure and let $TC$ and $FB$ be two collections of morphisms in $\mathcal C$ that satisfy Conditions \ref{2015.05.22.cond1}. Let further $Eq:(\wt{U};p)\sr U$ and $\Omega:\wt{U}\sr \wt{U}$ be a J1-structure and assume that the following conditions hold:
\begin{enumerate}
\item $U$ is fibrant,
\item $p$ is in $FB$,
\item $\omega$ is in $TC$.
\end{enumerate}
Then there exists a J-structure $Jp$ extending $(Eq,\Omega)$.
\end{theorem}
\begin{myproof}
Let us use the notations of Problem \ref{2015.05.12.l1}. We need to show that under the assumptions of the current theorem there exists a morphism that splits the square of Problem \ref{2015.05.12.l1} into two commutative triangles. Observe first that constructing such a splitting is equivalent to constructing the splitting of the square
$$
\begin{CD}
(\wt{U},p)\times_U (Fp,pFp) @>\sigma\circ adj(pr_2)>> U\times\wt{U}\\
@V\omega \times Id_{Fp} VV @VV Id_U\times p V\\
(E\wt{U},pE\wt{U})\times_U (Fp,pFp) @>\sigma'\circ adj(pr_1)>> U\times U
\end{CD}
$$
where
$$\sigma:(\wt{U},p)\times_U (Fp,pFp)\sr (Fp,pFp)\times_U  (\wt{U},p)$$
$$\sigma':(E\wt{U},pE\wt{U})\times_U (Fp,pFp) \sr  (Fp,pFp)  \times_U  (E\wt{U},pE\wt{U})$$
are permutations of the factors. 

It is easy to show that $U\times U$ is fibrant. Therefore it is sufficient to show that $Id_U\times p$ is in $FB$ and $\omega\times_U Id_{Fp}$ is in $TC$. The first fact follows from the assumption that $p$ is in $FB$ and that $U$ is fibrant. The obtain the second fact let us apply condition (3) on the classes $FB$ and $TC$ to $B=U$, $f=pE\wt{U}$, $i=\omega$ and $p=pFp$.  It remains to show that $pFp$ is in $FB$. We can represent $pFp$ as the composition
$$Fp\stackrel{pr_1}{\sr} I_{pE\wt{U}}(U) \stackrel{prI_{pE\wt{U}}}{\sr} U$$
The morphism $pE\wt{U}$ is in $FB$ as a composition of pull-backs of $p$ with respect to morphisms with fibrant domains through repeated application of Lemmas \ref{2015.05.14.l2} and \ref{2015.05.14.l4}. Therefore, the morphism $prI_{pE\wt{U}}$ is in $FB$ by Lemma \ref{2015.05.14.l1} and as a corollary we know that $I_{pE\wt{U}}(U)$ is fibrant. Similarly $I_p(U)$ is fibrant and $i_p(p)$ is in $FB$ and applying again Lemma \ref{2015.05.14.l2} we see that $pr_1$ is in $FB$. And again by Lemma \ref{2015.05.14.l4} we see that $pFp$ is in $FB$ which finishes the proof of the theorem.  
\end{myproof}
\begin{cor}
\llabel{2015.05.18.cor1}
Let $\cal C$ be a locally cartesian closed category with a Quillen model structure, $p$ a universe in $\cal C$ and $(Eq,\Omega)$ a J1-structure on $p$. Assume further that $p$ is a fibration and $\omega$ is a trivial cofibration and that in addition one of the following two conditions holds:
\begin{enumerate}
\item consider morphisms $f: B'\sr B$, $p_1:E_1\sr B$, $p_2:E_2\sr B$ and $i:E_1\sr E_2$ such that $p_1,p_2$ are fibrations and $i$ a trivial cofibration. Then the morphism
$$Id_{B'}\times i: (B',f)\times_B(E_1,p_1)\sr (B',f)\times_B(E_2,p_2)$$
is a trivial cofibration,
\item $U$ is fibrant and the pull-back of a trivial cofibration along a fibration is a trivial cofibration.
\end{enumerate}
Then $(Eq,\omega)$ can be extended to a full J-structure on $p$. 
\end{cor}

The following result can be used to produce many examples of universes with J-structures (but not the univalent universes). Let $\cal C$ be a locally cartesian closed category with coproducts of sequences $\amalg_{n\in\nn}X_n$. We let $in_n:X_n\sr \amalg_n X_n$ and $\langle f_n \rangle_{n} : \amalg_n X_n\sr Y$ denote the canonical morphisms. We let $\amalg f_n : \amalg_n X_n \amalg_n Y_n$ denote the morphism $\langle f_n\circ in_n \rangle_{n}$. 

Assume that these coproducts satisfy the following two conditions:
\begin{enumerate}
\item for a sequence of morphisms $f_n:E_n\sr B_n$ the square
$$
\begin{CD}
\amalg_n (E_n,f_n)\times_{B_n}(E_n,f_n) @>\amalg_n pr_2>> \amalg_n E_n\\
@V \amalg_n pr_1 VV @VV \amalg_n f_n V\\
\amalg_n E_n @>\amalg_n f_n >> \amalg_n B_n
\end{CD}
$$
is a pull-back square,
\item for a sequence of morphisms $f_n:E_n\sr B_n$ the square
$$
\begin{CD}
\amalg_n E_{n+1} @> \langle in_{n+1} \rangle_{n} >> \amalg_n E_n\\
@V\amalg_n f_{n+1} VV @VV \amalg_n f_n V\\
\amalg_n B_{n+1} @> \langle in_{n+1} \rangle_{n} >> \amalg_n B_n
\end{CD}
$$
is a pull-back square. 
\end{enumerate}
\begin{problem}
\llabel{2015.05.22.th2}
Let $\cal C$ be as above $FB$ and $TC$ two classes of morphisms satisfying one of the sets of conditions \ref{2015.05.22.cond1} or \ref{2015.05.22.cond2}. Assume in addition the following:
\begin{enumerate}
\item the coproduct of a sequence of morphisms from $TC$ is in $TC$ and the coproduct of a sequence of morphisms from $FB$ is in $FB$,
\item the composition of a morphism from $TC$ with an isomorphism is in $TC$,
\item  for any morphism $f: X \sr Y$ there is given an object $P(f)$ and morphisms $i_f:X\sr P(f)$, $q_f:P(f)\sr Y$ such that $i_f\in TC$, $q_f\in FB$ and $f=i_f\circ q_f$.
\end{enumerate}
To construct, for any universe $p: \wt{U}\sr U$ such that $p\in FB$ a sequence of morphisms $p_n:\wt{U}_n\sr U_n$ such that $p_0=p$, $p_n\in FB$ and $\amalg_n p$, with the universe structure defined by the fiber squares of $\cal C$, has a J-structure with $\omega\in TC$.
\end{problem}
\begin{construction}\rm\llabel{2015.05.23.constr1}
Define $p_n:\wt{U}_n\sr U_n$ inductively as follows. For $n=0$ we take $p_0=p$. To define $p_{n+1}$ consider the diagonal $\Delta_n:\wt{U}_n\sr (\wt{U}_n,p_n)\times_{U_n}(\wt{U}_n,p_n)$ and let 
$$p_{n+1}=q_{\Delta_n}:P(\Delta_n)\sr (\wt{U}_n,p_n)\times_{U_n}(\wt{U}_n,p_n)$$
so that in particular $U_{n+1}=(\wt{U}_n,p_n)\times_{U_n}(\wt{U}_n,p_n)$. 

Let $U_*=\amalg_n U_n$, $\wt{U}_*=\amalg_n\wt{U}_n$ and $p_*=\amalg_n p_n$. According to the first of the two properties that we required from the coproducts the canonical morphism 
$$\iota:\amalg_n (\wt{U}_n,p_n)\times_{U_n}(\wt{U}_n,p_n)\sr (\wt{U}_*,p_*)\times_{U_*}(\wt{U}_*,p_*)$$
is an isomorphism. Together with the second property applied to the right-most square this gives us a diagram with pull-back squares of the form:
$$
\begin{CD}
\amalg_n\wt{U}_{n+1} @>=>> \amalg_n\wt{U}_{n+1} @>=>> \amalg_n\wt{U}_{n+1} @> \langle in_{n+1} \rangle_{n} >> \wt{U}_*\\
@V r\circ \iota VV @V r VV @V r VV @VV p_* V\\
(\wt{U}_*,p_*)\times_{U_*}(\wt{U}_*,p_*) @>\iota^{-1}>> \amalg_n (\wt{U}_n,p_n)\times_{U_n}(\wt{U}_n,p_n) @>=>> \amalg_n U_{n+1} @> \langle in_{n+1} \rangle_{n} >> U_*
\end{CD}
$$
where $r=\amalg_n p_{n+1}$. Define $Eq$ as the composition of the lower horizontal arrows of this diagram (up to an isomorphism this is just $\langle in_{n+1}\rangle_{n}$). Since the squares of the diagram are pull-back the natural morphism
$$\iota':\amalg_n\wt{U}_{n+1}\sr ((\wt{U}_*,p_*)\times_{U_*}(\wt{U}_*,p_*),Eq)_{U_*} (\wt{U}_*,p_*)$$
is an isomorphism. Define
$$\Omega=(\amalg_n i_{\Delta_n})\circ \iota'\circ \langle in_{n+1} \rangle_{n}$$
such that then 
$$\omega=(\amalg_n i_{\Delta_n})\circ \iota'$$
By our assumptions $\omega\in TC$ and then by Theorem \ref{2015.05.22.th1} if $FB$ and $TC$ satisfied Conditions \ref{2015.05.22.cond2} or by Theorem \ref{2015.05.16.th1} if $FB$ and $TC$ satisfied Conditions \ref{2015.05.22.cond1} we conclude that $(Eq,\Omega)$ can be extended to a full J-structure on $p_*$.
\end{construction}

\subsection{Constructing a J-structure on $CC({\mathcal C},p)$ from a J-structure on $p$}
The construction of a C-system $CC({\mathcal C},p)$ from a category with a universe $p$ and a final object $pt$ was presented in \cite{Cfromauniverse} and summarized in \cite{fromunivwithPi}. Let us recall some facts and notations. The underlying category of $CC({\mathcal C},p)$ is equipped with a functor $int$ to $\mathcal C$. Note that while $int$ is the identity on morphisms by construction of $CC({\mathcal C},p)$, the notations for the same element of $Hom(\Gamma',\Gamma)$ and $Hom(int(\Gamma'),int(\Gamma))$ may differ. In particular for $f:\Gamma'\sr \Gamma$ and $F:\Gamma\sr U$ we have
\begin{eq}
\llabel{2015.04.02.eq2}
q(f,int(\Gamma,F))=Q(f,F)
\end{eq}

For each $\Gamma\in Ob(CC({\mathcal C},p))$ we have natural bijections 
\begin{eq}
\llabel{2015.03.27.eq7b}
u_1:Ob_1(\Gamma)\sr Hom(int(\Gamma),U)
\end{eq}
\begin{eq}
\llabel{2015.03.27.eq7a}
\wt{u}_1:\wt{Ob}_1(\Gamma)\sr Hom(int(\Gamma),\wt{U})
\end{eq}
where $u_1^{-1}(F)=(\Gamma,F)$ and where
\begin{eq}
\llabel{2015.03.31.eq5}
\wt{u}_1(s)=s\circ Q(u_1(\partial(s)))
\end{eq}
In particular,
$$\wt{u}_1(s)\circ p=s\circ Q(u_1(\partial(s)))\circ p=s\circ p_{\partial(s)}\circ u_1(\partial(s))=u_1(\partial(s))$$
i.e., with respect to these bijections the function $\partial:\wt{Ob}_1(\Gamma)\sr Ob_1(\Gamma)$ is given by composition with $p:\wt{U}\sr U$. 
\begin{problem}
\llabel{2015.03.27.prob3}
Let $Eq:(\wt{U};p)\sr U$ be a J0-structure on a universe $p$ in a category $\mathcal C$. To construct a J0-structure on $CC({\mathcal C},p)$.
\end{problem}
\begin{construction}\rm
\llabel{2015.03.27.constr3}
Since the canonical squares are pull-back squares bijections $u_1$ and $\wt{u}_1$ gives us a bijection
$$\wt{uu}:\{o,o'\in\wt{Ob}_1(\Gamma)\,|\,\partial(o)=\partial(o')\} \sr Hom(int(\Gamma),(\wt{U};p))$$
where $\wt{uu}(o,o')=\wt{u}_1(o)*\wt{u}_1(o')$. We set:
$$IdT(o,o')=u_1^{-1}(\wt{uu}(o,o')\circ Eq).$$
\end{construction}
We let $IdT_{Eq}$ denote the J0-structure on $CC({\mathcal C},p)$ constructed from $Eq$ in Construction \ref{2015.03.27.constr3}. Note that
\begin{eq}
\llabel{2015.03.31.eq1}
int(IdT(o,o'))=(int(\Gamma);(\wt{u}_1(o)*\wt{u}_1(o'))\circ Eq)
\end{eq}
Recall that in \cite{Csubsystems} we let $p_{\Gamma,n}:\Gamma\sr ft^n(\Gamma)$ denote the composition of $n$ canonical projections $p_{\Gamma}\circ \dots\circ p_{ft^{n-1}(\Gamma)}$. 
\begin{lemma}
\llabel{2015.03.27.l1}
Let $Eq$ be a J0-structure on $p$. Let $\Gamma\in Ob$ and $F:int(\Gamma)\sr U$. Then one has:
$$int(IdxT(\Gamma,F))=(int(\Gamma);F)_{E}$$
$$p_{IdxT(T),3} = p^E_{\Gamma,F}$$
$$Q(F)_{E}\circ Q(Eq)=Q(Q(Q(F)\circ p)\circ Eq)$$
\end{lemma}
\begin{myproof}
Let $T=(\Gamma,F)$. We have:
$$int(IdxT(T))=int(IdT_{p_T^*(T)}(o,o'))=(int(p_T^*(T));(\wt{u}_1(o)*\wt{u}_1(o'))\circ Eq)$$
where 
$$o=p_{p_T^*(T)}^*(\delta(T))$$
$$o'=\delta(p_T^*(T))$$
We have
$$int(p_T^*(T))=(int(\Gamma);F,Q(F)\circ p)$$
and
$$\wt{u}_1(p_{p_T^*(T)}^*(\delta(T)))=p_{(int(\Gamma,F, Q(F)\circ p)}\circ Q(F)$$
$$\wt{u}_1(\delta(p_T^*(T)))=Q(Q(F)\circ p)$$
which shows that $\wt{u}_1(o)*\wt{u}_1(o')=Q(Q(F),p)$ and completes the proof of the first and the second equations.

The third equality is a corollary of the equality $Q(F)_{E}=Q(Q(Q(F),p),Eq)$ and the equality $Q(f,F)\circ Q(F)=Q(f\circ F)$. 
\end{myproof}
\begin{problem}
\llabel{2015.03.27.prob4}
Let $Eq:(\wt{U};p)\sr U$, $\Omega:\wt{U}\sr \wt{U}$ be a J0-structure and a J1-structures on a universe $p$ in a category $\mathcal C$. To construct a J1-structure $refl(Eq,\Omega)$ over $IdT_{Eq}$ on $CC({\mathcal C},p)$.
\end{problem}
\begin{construction}\rm
\llabel{2015.03.27.constr4}
Due to the natural bijections  (\ref{2015.03.27.eq7a}) the morphism $\Omega$ defines maps
$$refl:\wt{Ob}_1(\Gamma)\sr \wt{Ob}_1(\Gamma)$$
by the formula
$$refl(s)=\wt{u}_1^{-1}(\wt{u}_1(s)\circ \Omega)$$
that are natural in $\Gamma$. The equation (\ref{2015.03.27.eq8}) follows from the commutativity of the square (\ref{2015.03.27.sq1}).
\end{construction}
We let $refl_{\Omega}$ denote the J1-structure constructed from $\Omega$ in Construction \ref{2015.03.27.constr4}. 
The following technical lemma is only needed in the proof of Lemma \ref{2015.03.31.l2}.
\begin{lemma}
\llabel{2015.04.02.l3}
For $s\in \wt{Ob}_1(\Gamma)$ one has:
$$refl_{\Omega}(s)\circ Q(s\circ Q(F)\circ \Omega\circ p)=s\circ Q(F)\circ \Omega$$
where $F=u_1(\partial(s))$.
\end{lemma}
\begin{myproof}
We have
$$u_1(\partial(refl_{\Omega}(s)))=\wt{u}_1(refl(s))\circ p=\wt{u}_1(s)\circ \Omega\circ p=s\circ Q(F)\circ \Omega\circ p$$
therefore
$$\wt{u}_1(refl_{\Omega}(s))=refl_{\Omega}(s)\circ Q(u_1(\partial(refl_{\Omega}(s))))=$$
$$refl_{\Omega}(s)\circ Q(s\circ Q(F)\circ \Omega\circ p)$$
On the other hand, by definition of $refl_{\Omega}$,
$$\wt{u}_1(refl_{\Omega}(s))=\wt{u}_1(s)\circ \Omega=s\circ Q(F)\circ \Omega.$$
\end{myproof}

\begin{lemma}
\llabel{2015.03.31.l2}
Given $Eq$ and $\Omega$ consider the corresponding $IdT$ and $refl$. For $T\in Ob_1(\Gamma)$ let %
$$rf_T:int(T)\sr int(IdxT(T))$$
be the morphism constructed in Construction \ref{2015.03.27.constr2}. On the other hand let 
$$F^*(\omega):(int(\Gamma);F)\sr (int(\Gamma);F)_{E}$$
is the pull-back of $\omega:\wt{U}\sr E\wt{U}$ with respect to $F=u_1(T)$ i.e. the unique morphism 
$$(int(\Gamma);F)\sr (int(\Gamma);F)_{E}$$
such that 
$$F^*(\omega)\circ p^{E}_{int(\Gamma),F}=p_{T}$$
$$F^*(\omega)\circ Q(int(\Gamma),F)_{E}=Q(F)\circ \omega$$
Then 
$$rf_T=F^*(\omega)$$
\end{lemma}
\begin{myproof}
In view of Lemma \ref{2015.03.27.l1}, both $rf_T$ and $F^*(\omega)$ are morphisms from $(int(\Gamma);F)$ to $(int(\Gamma);F)_{E}$. Let us denote $int(\Gamma)$ by $X$ and $(int(\Gamma);F,Q(F)\circ p)$ by $Y$. We have
$$(X;F)_{E}=(Y;Q(Q(F),p)\circ Eq)$$
and we can see this object as a part of the diagram with two pull-back squares:
$$
\begin{CD}
(Y;Q(Q(F),p)\circ Eq) @>h_1>> E\wt{U} @>h_2>> \wt{U}\\
@Vp_{Y, Q(Q(F),p)\circ Eq}VV @VVV @VVpV\\
Y @>Q(Q(F),p)>> (\wt{U},p) @>Eq>> U
\end{CD}
$$
We have two projections
$$h=h_1\circ h_2=Q(Q(Q(F),p)\circ Eq):(X;F)_{E}\sr \wt{U}$$
$$v:p_{Y,Q(Q(F),p)\circ Eq}:(X;F)_{E} \sr Y$$
We need to check that
$$rf_T\circ h=F^*(\omega)\circ h$$
$$rf_T\circ v =F^*(\omega)\circ v$$
The morphism $rf_T$ is defined in (\ref{2015.04.02.eq1}) as 
$$rf_T=refl(\delta(T))\circ q(\delta(T),IdxT(T))=refl(\delta(T))\circ Q(\delta(T),Q(Q(F),p)\circ Eq)$$
where the second equation is from (\ref{2015.04.02.eq2}).
We have 
$$rf_T\circ h = refl(\delta(T))\circ Q(\delta(T),Q(Q(F),p)\circ Eq)\circ h=refl(\delta(T))\circ Q(\delta(T)\circ Q(Q(F),p)\circ Eq)=$$$$refl(\delta(T))\circ Q(Q(F)\circ \Delta\circ Eq)=refl(\delta(T))\circ Q(Q(F)\circ (\omega\circ pE\wt{U})\circ Eq)=$$
$$refl(\delta(T))\circ Q(Q(F)\circ \omega\circ Q(Eq)\circ p)=refl(\delta(T))\circ Q(Q(F)\circ \Omega\circ p)$$
On the other hand
$$F^*(\omega)\circ h=F^*(\omega)\circ h_1\circ h_2=Q(F)\circ \omega\circ h_2=Q(F)\circ \omega\circ Q(Eq)=Q(F)\circ \Omega$$
We have $u_1(\partial(\delta(T)))=Q(F)\circ p$ and 
$$\delta(T)\circ Q(u_1(\partial(\delta(T))))=\delta(T)\circ Q(Q(F)\circ p)= \delta(T)\circ Q(p_T\circ F)=$$
$$\delta(T)\circ q(p_T,T)\circ Q(F)=Id_{int(T)}\circ Q(F)=Q(F)$$
Therefore by Lemma \ref{2015.04.02.l3} we have
$$refl(\delta(T))\circ Q(Q(F)\circ \Omega\circ p)=Q(F)\circ \Omega$$
Which proves that $rf_T\circ h=F^*(\omega)\circ h$. 

We have $rf_T\circ v=\delta(T)$ because the square (\ref{2015.03.31.eq3}) commutes.  Both $rf_T\circ v$ and $F^*(\omega)\circ v$ are morphisms $int(T)\sr int(p_T^*(T))$. Since $p_T^*(T)$ is a part of a pull-back square with the projections being $p_{p_T^*(T)}$ and $Q(Q(F)\circ p)$ we need to check that
$$\delta(T)\circ Q(Q(F)\circ p)=  F^*(\omega)\circ v\circ Q(Q(F)\circ p)=F^*(\omega)\circ h_1\circ pE\wt{U}=$$
$$Q(F)\circ \omega \circ pE\wt{U}=Q(F)\circ \Delta$$
which holds by a simple computation, and
$$Id_{int(T)}=\delta(T)\circ p_{p_T^*(T)}=F^*(\omega)\circ v\circ p_{p_T^*(T)}$$
For this equality we need to verify two further equalities
$$Q(F)=Id_{int(T)}\circ Q(F)=F^*(\omega)\circ v\circ p_{p_T^*(T)}\circ Q(F)$$
and
$$p_T=Id_{int(T)}\circ p_T=F^*(\omega)\circ v\circ p_{p_T^*(T)}\circ p_T$$
The second one is the second equality of the two that define $F^*(\omega)$. For the first one we have
$$F^*(\omega)\circ v\circ (p_{p_T^*(T)}\circ Q(F))=F^*(\omega)\circ v\circ Q(Q(F),p)\circ p_{\wt{U},p}=$$
$$Q(F)\circ \Delta\circ p_{\wt{U},p}=Q(F).$$
This completes the proof of Lemma \ref{2015.03.31.l2}.
\end{myproof}

\comment{ 
\begin{problem}
\llabel{2015.05.04.prob2}
Let $Eq$ be a J0-structure on a universe $p$ in $\cal C$. To construct, for all objects $\Gamma\in Ob(CC({\mathcal C},p))$ two maps
$$\psi:D_{pE\wt{U}}(int(\Gamma),U)\sr Ob(CC({\mathcal C},p))$$
$$\wt{\psi}:D_{pE\wt{U}}(int(\Gamma),\wt{U})\sr \wt{Ob}(CC({\mathcal C},p))$$
such that 
\begin{eq}
\llabel{2015.05.04.eq3}
\partial(\wt{\psi}(d))=\psi(D_{pE\wt{U}}(int(\Gamma),p)(d))
\end{eq}
%
\end{problem}
\begin{construction}
\llabel{2015.05.04.constr2}\rm
Let us construct $\wt{\psi}$.  An element $(F,\wt{G})$ of $D_{pE\wt{U}}(int(\Gamma),\wt{U})$ is a pair where $F:int(\Gamma)\sr U$ and $\wt{G}:(int(\Gamma);F)_{E}\sr \wt{U}$. By Lemma \ref{2015.03.27.l1} we have $(int(\Gamma),F)_{E}=int(IdxT(\Gamma,F))$. Therefore $\wt{u}_{1,IdxT(\Gamma,F)}^{-1}(\wt{G})$ is defined and belongs to 
$\wt{Ob}_1(IdxT(\Gamma,F))$. We set 
$$\wt{\psi}(F,\wt{G})=\wt{u}_{1,IdxT(\Gamma,F)}^{-1}(\wt{G})$$
Similarly we set 
$$\psi(F,G)=u_{1,IdxT(\Gamma,F)}^{-1}(G)$$
The proof of the relation (\ref{2015.05.04.eq3}) is straightforward. 
\end{construction}
}
\begin{problem}
\llabel{2015.04.04.prob1}
Let $(Eq,\Omega,Jp)$ be a J-structure on a universe $p$. To construct for all $\Gamma\in Ob=Ob(CC({\mathcal C},p))$, for all $T\in Ob_1(\Gamma)$, for all $P\in Ob_1(IdxT(T))$, for all $s0\in \wt{Ob}(rf_T^*(P))$, an element $J(\Gamma,T,P,s0)$ of $\wt{Ob}(P)$.
\end{problem} 
\begin{construction}\rm
\llabel{2015.04.04.constr1}
Let $X:=int(\Gamma)$, $F:=u_1(T):X\sr U$. By Lemma \ref{2015.03.27.l1} we have $int(IdxT(T))=int(IdxT(\Gamma,F))=(int(\Gamma);F)_{E}$. Therefore we further have
$G:=u_1(P):(int(\Gamma);F)_{E}\sr U$ and $\wt{H}:=\wt{u}_1(s0):(X;F)\sr \wt{U}$. 

Let us show first that
$$\eta_{pE\wt{U}}(F,G)\circ I^{\omega}(U)=\eta_p(F,\wt{H})\circ I_p(p)$$
Indeed,
$$\eta_{pE\wt{U}}(F,G)\circ I^{\omega}(U)=\eta_p(F,G\circ F^*(\omega))=\eta_p(F,rf_T\circ G)=$$
$$\eta_p(F,s0\circ Q(rf_T\circ G)\circ p)=\eta_p(F,\wt{H}\circ p)=\eta_p(F,\wt{H})\circ I_p(p)$$
where the first equality is from Lemma \ref{2015.04.02.l4}, the second from Lemma \ref{2015.03.31.l2}, the third from the commutativity of the canonical square and the fact that $s0$ is a section, the fourth from (\ref{2015.03.31.eq5}) and the fifth from naturality of $\eta_{p,X,V}$ in $V$. 

Therefore the pair $(\eta_{pE\wt{U}}(F,G),\eta_p(F,\wt{H}))$ gives us a morphism
$$\phi(\Gamma,T,P,s0):X\sr (I_{pE\wt{U}}(U), I^{\omega}(U)) \times_{I_p(U)} (I_p(\wt{U}),I_p(p))$$
and compositing it with $Jp$ (cf. Definition \ref{2015.03.27.def6}) we obtain a morphism 
$$\phi(\Gamma,T,P,s0)\circ Jp: X\sr I_{pE\wt{U}}(\wt{U})$$
Consider the element 
$$(F_1,F_2)=\eta^!_{pE\wt{U}}(\phi(\Gamma,T,P,s0)\circ Jp)\in D_{pE\wt{U}}(X,\wt{U})$$
By \cite[Problem 3.8(1)]{fromunivwithPi} we have 
$$(F_1,F_2\circ p)=D_{pE\wt{U}}(X,p)(F_1,F_2)=\eta^!_{pE\wt{U}}(\phi(\Gamma,T,P,s0)\circ Jp\circ I_{pE\wt{U}}(p))=$$
$$\eta^!_{pE\wt{U}}(\eta_{pE\wt{U}}(F,G))=(F,G)$$
Therefore, $F_2$ is of the form $(X;F)_E\sr \wt{U}$ i.e. of the form $\wt{u}_{1,IdxT(T)}(J))$ for some $J$ such that $\partial(J)=P$. 
\end{construction}
\begin{remark}\rm
\llabel{2015.05.08.rem1}
Note that the defining property of $J=J_{Jp}(\Gamma,T,P,s0)$ is that it is the unique element of $\wt{Ob}(CC({\mathcal C},p))$ that satisfies the equation
$$\eta_{pE\wt{U}}(u_{1,\Gamma}(T),\wt{u}_{1,IdxT(T)}(J))=\phi(\Gamma,T,P,s0)\circ Jp$$
where 
$$\phi(\Gamma,T,P,s0):int(\Gamma)\sr (I_{pE\wt{U}}(U),I^{\omega}(U))\times_{I_p(U)}(I_p(\wt{U}),I_p(p))$$
is given by the pair of morphisms $(\eta_{pE\wt{U}}(u_{1,\Gamma}(T), u_{1,IdxT(T)}(P)), \eta_p(u_{1,\Gamma}(T),\wt{u}_{1,\Gamma}(s0)))$.
\end{remark}

\begin{lemma}
\llabel{2015.04.04.l4}
Let $Eq$ be a J0-structure on a universe $p$, $f:\Gamma\sr \Gamma'$ a morphism in $CC({\mathcal C},p)$ and $F:int(\Gamma)\sr U$ a morphism in $\mathcal C$.  Let $q3:int(IdxT(\Gamma',f\circ F))\sr int(IdxT(\Gamma,F))$ be the morphism $q(f,IdxT(\Gamma,F),3)$ defined by $\Gamma$ since $ft^3(IdxT(\Gamma,F))=\Gamma$. Then $q3=Q(f,F)_{E}$.
\end{lemma}
\begin{myproof}
Let $X:=int(\Gamma)$ and $X':=int(\Gamma')$. By definition, $Q(f,F)_{E}$ is the  unique morphism such that
$$Q(f,F)_{E}\circ Q(F)_{E}=Q(f\circ F)_{E}$$
$$Q(f,F)_{E}\circ p^{E}_{X,F}=p^{E}_{X',f\circ F}\circ f$$
We will be building the proof using the following diagram
$$
\begin{CD}
(X',f\circ F)_{E} @>Q(Q(Q(f,F),Q(F)\circ p),Q(Q(F),p)\circ Eq)>> (X;F)_{E} @>Q(Q(Q(F),p),Eq)>> E\wt{U} @>Q(Eq)>> \wt{U}\\
@Vp_3VV @VVV @VVp_1V @VVp V\\
\BB @>Q(Q(f,F),Q(F)\circ p)>> \BB @>Q(Q(F),p)>> \BB @>Eq>>U\\
@| @| @|\\
\BB @>Q(Q(f,F),Q(F)\circ p)>> \BB @>Q(Q(F),p)>> \BB @>Q(p)>> \wt{U}\\
@VVV @VVV @VVp_2V @VVp V\\
\BB @>Q(f,F)>> \BB @>Q(F)>> \BB @>p>> U\\
@| @| @|\\
\BB @>Q(f,F)>> \BB @>Q(F)>> \wt{U}\\
@VVV @VVV @VVp V\\
X' @>f>> X @>F>> U
\end{CD}
$$
By construction that is seen on this diagram we have:
$$q3=Q(Q(Q(f,F),Q(F)\circ p),Q(Q(F),p)\circ Eq)$$
$$Q(X,F)_{E}=Q(Q(Q(F),p),Eq)$$
and
$$Q(X', f\circ F)_{E}=Q(Q(Q(f\circ F),p),Eq)$$
Therefore, the first equation that we need to verify is
$$Q(Q(Q(f,F),Q(F)\circ p),Q(Q(F),p)\circ Eq)\circ Q(Q(Q(F),p),Eq)=Q(Q(Q(f\circ F),p),Eq)$$
By \cite[Lemma 3.2]{fromunivwithPi} we have, together with the defining rule $Q(a,A)\circ Q(A)=Q(a\circ A)$, also the rule:
$$Q(a_1,a_2\circ A)\circ Q(a_2,A)=Q(a_1\circ a_2, A)$$
Applying it twice and then the defining rule we get:
$$Q(Q(Q(f,F),Q(F)\circ p),Q(Q(F),p)\circ Eq)\circ Q(Q(Q(F),p),Eq)=$$
$$Q(Q(Q(f,F),Q(F)\circ p)\circ Q(Q(F),p), Eq)=$$
$$Q(Q(Q(f,F)\circ Q(F),p),Eq)=Q(Q(Q(f\circ F),p),Eq)$$
which gives us the first equation. The second equation is immediate from the commutativity of the three squares that define $q3$. 
\end{myproof}
\begin{lemma}
\llabel{2015.04.04.l1}
Let $(Eq,\Omega,Jp)$ be a J-structure on a universe $p$. Then the morphisms of Construction \ref{2015.04.04.constr1} are natural in $\Gamma$ i.e. for any $f:\Gamma'\sr \Gamma$ one has
\begin{eq}\llabel{2015.04.04.eq3}
f^*(J_{Jp}(\Gamma,T,P,s0))=J_{Jp}(\Gamma',f^*(T),f^*(P),f^*(s0))
\end{eq}
\end{lemma}
\begin{myproof}
Let us write $J$ for $J_{Jp}(\Gamma,T,P,s0)$ and $J'$ for $J_{Jp}(\Gamma',f^*(T),f^*(P),f^*(s0))$ and use the notations of Construction \ref{2015.04.04.constr1}. Recall that for $f:\Gamma'\sr \Gamma$ the operation $f^*$ is defined only on $Ob_1(\Gamma)$. In all other uses it is an abbreviation for operations such as $X\mapsto f^*(X,i)$ and $s\mapsto f^*(s,i)$ for various $i$ (see \cite{Csubsystems}). In particular, (\ref{2015.04.04.eq3}) is an abbreviation for 
$$f^*(J(\Gamma,T,P,s0),4)=J(\Gamma',f^*(T),f^*(P,4),f^*(s0,2))$$
which in its turn translates into the equation in $\wt{Ob}_1(IdxT(f^*(T)))$ of the form
$$q(f,IdxT(T),3)^*(J,1)=J'$$
We have:
$$\eta_{pE\wt{U}}(F,\wt{u}_1(J))=\phi(\Gamma,T,P,s0)\circ Jp$$
$$\eta_{pE\wt{U}}(f\circ F, \wt{u}_1(J'))=\phi(\Gamma',f^*(T),f^*(P),f^*(s0))\circ Jp$$
By naturality of $\eta$ with respect to the first argument we have
$$f\circ \eta_{pE\wt{U}}(F,\wt{u}_1(J))=\eta_{pE\wt{U}}(f\circ F, Q(f,F)_{E}\circ \wt{u}_1(J))$$
Therefore, by Lemma \ref{2015.04.04.l4} we have
$$f\circ \eta_{pE\wt{U}}(F,\wt{u}_1(J))=\eta_{pE\wt{U}}(f\circ F, \wt{u}_1(Q(f,F)_{E}^*(J,1)))=$$
$$\eta_{pE\wt{U}}(f\circ F, \wt{u}_1(q(f,IdxT(T),3)^*(J,1)))$$
Since both $\eta_{pE\wt{U}}$ and $\wt{u}_1$ are bijections and in particular injections it is sufficient to show that
$$f\circ \phi(\Gamma,T,P,s0)\circ Jp = \phi(\Gamma',f^*(T),f^*(P),f^*(s0))\circ Jp$$
or that
$$f\circ \phi(\Gamma,T,P,s0)=\phi(\Gamma',f^*(T),f^*(P),f^*(s0))$$
Since both $\phi$ expressions are morphism into a product this amounts to two equations that, taking into account the definition of $\phi$ in Construction \ref{2015.04.04.constr1} are:
$$f\circ \eta_{pE\wt{U}}(F,G)=\eta_{pE\wt{U}}(f\circ F, u_1(f^*(P)))$$
and
$$f\circ \eta_p(F,\wt{H})=\eta_p(f\circ F, \wt{u}_1(f^*(s0)))$$
The first equality follows naturality of $\eta$ and Lemma \ref{2015.04.04.l4}. The second equality follows from naturality of $\eta$. This finished the proof of Lemma \ref{2015.04.04.l1}.
\end{myproof}
\begin{lemma}
\llabel{2015.04.04.l5}
Let $(Eq,\Omega,Jp)$ be a J-structure on a universe $p$. Then the morphisms of  Construction \ref{2015.04.04.constr1} satisfy the second condition of the definition of a J2-structure, i.e., for all $\Gamma$, $T$, $P$ and $s0$ as above one has
$$rf_T^*(J_{Jp}(\Gamma,T,P,s0))=s0$$
\end{lemma}
\begin{myproof}
Let $J=J_{Jp}(\Gamma,T,P,s0)$. Then, using the notations of Construction \ref{2015.04.04.constr1} we have:
$$\eta_{E\wt{U}}(F,\wt{u}_1(J))=\phi\circ Jp$$
Then
$$\eta_{E\wt{U}}(F,\wt{u}_1(J))\circ I^{\omega}(\wt{U})=\eta_p(F,F^*(\omega)\circ \wt{u}_1(J))$$
By Lemma \ref{2015.03.31.l2} we have $F^*(\omega)=rf_T$. Therefore,
$$\eta_{E\wt{U}}(F,\wt{u}_1(J))\circ I^{\omega}(\wt{U})=\eta_p(F,rf_T\circ \wt{u}_1(J))=\eta_p(F,\wt{u}_1(rf_T^*(J)))$$
On the other hand
$$\phi\circ Jp\circ I^{\omega}(\wt{U})=\phi\circ prI_p(\wt{U})$$
by (\ref{2015.04.04.eq1}) which equals, by construction, $\eta_p(F,\wt{u}_1(s0))$. Therefore,
$$\eta_p(F,\wt{u}_1(rf_T^*(J)))=\eta_p(F,\wt{u}_1(s0))$$
and using again that both $\eta$ and $\wt{u}_1$ are injective we conclude that $rf_T^*(J)=s0$. 
\end{myproof}
\begin{problem}
\llabel{2015.04.04.prob2}
Let $(Eq,\Omega,Jp)$ be a J-structure on a universe $p$. To construct a J-structure on $CC({\mathcal C},p)$ relative to $IdT_{Eq}$ and $refl_{\Omega}$.
\end{problem} 
\begin{construction}\rm
\llabel{2015.04.04.constr2}
One has to combine Construction \ref{2015.04.04.constr1} with Lemmas \ref{2015.04.04.l1} and \ref{2015.04.04.l5}.
\end{construction}

\section{Functoriality of the J-structures}

\subsection{A theorem about functors between categories with two universes}
\label{twouniv}
Before we can formulate the definition of what it means for a universe category functor to be compatible with J-structures we need some general results about functors between categories with two universes that we will later apply to the universes $p:\wt{U}\sr U$ and $pE\wt{U}:E\wt{U}\sr U$ in a locally cartesian closed category $\cal C$. 

Given two universes $(p,p_{X,F},Q(F))$ and $(p',p_{X,F}',Q(F)')$ where $p:\wt{U}\sr U$ and $p':\wt{U}'\sr U$ and the canonical squares are of the form
$$
\begin{CD}
(X;F) @>Q(F)>> \wt{U}\\
@Vp_{X,F} VV @VVp V\\
X @>F>> U
\end{CD}
\spc
\begin{CD}
(X;F)' @>Q(F)'>> \wt{U}'\\
@Vp_{X,F}' VV @VVp' V\\
X @>F>> U
\end{CD}
$$
and $f:\wt{U}'\sr \wt{U}$ over $U$, we let $F^*(f)$ denote the unique morphism $(X;F)'\sr (X;F)$ such that 
\begin{eq}\llabel{2015.04.08.eq3}
F^*(f)\circ Q(F)=Q(F)'\circ f
\end{eq}
\begin{eq}\llabel{2015.04.08.eq4}
F^*(f)\circ p_{X,F}=p_{X,F}'
\end{eq}
Note that $F^*(f)$ depends on the universe structures on $p$ and $p'$. Even when two universe structures give the same choices of the objects $(X;F)$ and $(X;F)'$ the difference in the choice of some of the morphisms, e.g., $Q(F)$ will affect morphisms $F^*(f)$. We will need the following lemma about these morphisms. 

For $X'\stackrel{f}{\sr}X \stackrel{F}{\sr}U$ we let $Q(f,F)$ denote the morphism
$$(p_{X',f\circ F}\circ f)*Q(f\circ F):(X';f\circ F)\sr (X;F)$$
We let $Q'(-)$ and $Q'(-,-)$ denote the morphisms $Q(-)$ and $Q(-,-)$ relative to the universe $p'$.
\begin{lemma}
\llabel{2015.04.20.l1}
Let $X'\stackrel{g}{\sr}X\stackrel{F}{\sr}U$ be two morphisms. Then the square
$$
\begin{CD}
(X';g\circ F)' @>Q'(g,F)>> (X;F)'\\
@V(g\circ F)^*(f)VV @VVF^*(f) V\\
(X';g\circ F) @>Q(g,F)>> (X;F)
\end{CD}
$$
commutes.
\end{lemma}
\begin{myproof}
Since $(X;F)$ is a fiber product relative to the projections $p_{X,F}$ and $Q(F)$ it is sufficient to verify that
$$Q'(g,F)\circ F^*(f)\circ Q(F)=(g\circ F)^*(f) \circ Q(g,F)\circ Q(F)$$
and
$$Q'(g,F)\circ F^*(f)\circ p_{X,F}=(g\circ F)^*(f) \circ Q(g,F)\circ p_{X,F}$$
which easily follows from the defining equations for $Q(-,-)$ and $(-)^*$.
\end{myproof}

Let $({\mathcal C},p,pt)$, $({\mathcal C}',p',pt')$ be two universe categories such that $\mathcal C$ and $\mathcal C'$ are equipped with locally cartesian closed structures. In \cite[Construction 5.6]{fromunivwithPi} we have defined, for any universe category functor
$${\bf\Phi}=(\Phi,\phi,\wt{\phi}):({\mathcal C},p,pt)\sr ({\mathcal C}',p',pt')$$
and any $V\in {\mathcal C}$,  a morphism
$$\chi_{\bf\Phi}(V):\Phi(I_p(V))\sr I_{p'}(\Phi(V))$$
We now need to consider the case when we have the following collection of data:
\begin{enumerate}
\item two universes $p_1$, $p_2$ in $\mathcal C$ with the common codomain $U$ and a morphism $g:\wt{U}_1\sr \wt{U}_2$ over $U$,
\item two universes $p_1'$, $p_2'$ in $\mathcal C'$ with the common codomain $U'$ and a morphism $g':\wt{U}_1'\sr \wt{U}_2'$ over $U'$,
\item a functor $\Phi:{\mathcal C}\sr {\mathcal C}'$,
\item a morphism $\phi:\Phi(U)\sr U'$,
\item two morphisms $\wt{\phi}_i:\Phi(\wt{U}_i)\sr \wt{U}_i'$, $i=1,2$
\end{enumerate}
and this data is such that:
\begin{enumerate}
\item the square
$$
\begin{CD}
\Phi(\wt{U}_1) @>\wt{\phi}_1>> \wt{U}_1'\\
@V\Phi(g) VV @VV g'V\\
\Phi(\wt{U}_2) @>\wt{\phi}_2>> \wt{U}_2'
\end{CD}
$$ 
commutes,
\item the triples ${\bf\Phi}_i:=(\Phi,\phi,\wt{\phi}_i)$, $i=1,2$ are universe category functors i.e. 
$\Phi$ takes canonical squares of $p_1$ and $p_2$ to pull-back squares and the squares
$$
\begin{CD}
\Phi(\wt{U}_1) @>\wt{\phi}_1>> \wt{U}_1'\\
@V\Phi(p_1) VV @VVp_1' V\\
\Phi(U) @>>\phi> U
\end{CD}
\spc\spc\spc
\begin{CD}
\Phi(\wt{U}_2) @>\wt{\phi}_2>> \wt{U}_2'\\
@V\Phi(p_2) VV @VVp_2' V\\
\Phi(U) @>>\phi> U
\end{CD}
$$
are pull-back squares. 
\end{enumerate}
Let us denote the exchange morphisms
$$\chi_{\bf\Phi_i}(V):\Phi(I_{p_i}(V))\sr I_{p_i}(\Phi(V))$$
by $\chi_i(V)$. The maps ${\bf\Phi}_i^2$ in the following lemma are constructed in \cite[Construction 5.2]{fromunivwithPi}. 
\begin{lemma}
\llabel{2015.04.08.l1}
Under the previous assumptions and notations the squares
$$
\begin{CD}
D_{p_2}(X,V) @>{\bf\Phi}_2^2>> D_{p_2}(\Phi(X),\Phi(V))\\
@V{D^g(X,V)} VV @VV {D^{g'}(\Phi(X),\Phi(V))}V\\
D_{p_1}(X,V) @>{\bf\Phi}_1^2>> D_{p_1'}(\Phi(X),\Phi(V))
\end{CD}
$$
commute.
\end{lemma}
\begin{myproof}
Let $(F_1,F_2)\in D_{p_2}(X,F)$ then
$$D^{g'}(\Phi(X),\Phi(V))({\bf\Phi}_2^2(F_1,F_2))=D^{g'}(\Phi(X),\Phi(V))(\Phi(F_1)\circ \phi, \iota_2\circ \Phi(F_2))=$$
$$=(\Phi(F_1)\circ \phi, (\Phi(F_1)\circ \phi)^*(g')\circ \iota_2\circ \Phi(F_2))$$
On the other hand
$${\bf\Phi}_1^2(D^g(X,V)(F_1,F_2))={\bf\Phi}_1^2(F_1,F_1^*(g)\circ F_2)=$$
$$(\Phi(F_1)\circ\phi,\iota_1\circ \Phi(F_1^*(g)\circ F_2))$$
It remains to check that 
$$(\Phi(F_1)\circ \phi)^*(g')\circ \iota_2\circ \Phi(F_2)=\iota_1\circ \Phi(F_1^*(g)\circ F_2)$$
For which it is sufficient to check that
$$(\Phi(F_1)\circ \phi)^*(g')\circ \iota_2=\iota_1\circ \Phi(F_1^*(g))$$
The codomain of both morphisms is $\Phi((X;F_1)_{p_2})$ and since $\Phi$ takes canonical squares based on $p_2$ to pull-back squares it is sufficient to check that 
$$(\Phi(F_1)\circ \phi)^*(g')\circ \iota_2\circ \Phi(Q_2(F_1))\circ\wt{\phi}_2=\iota_1\circ \Phi(F_1^*(g))\circ \Phi(Q_2(F_1))\circ\wt{\phi}_2$$
and
$$(\Phi(F_1)\circ \phi)^*(g')\circ \iota_2\circ \Phi(p_{X,F_1,2})=\iota_1\circ \Phi(F_1^*(g))\circ \Phi(p_{X,F_1,2})$$
For the first equation we have
$$(\Phi(F_1)\circ \phi)^*(g')\circ \iota_2\circ \Phi(Q_2(F_1))\circ\wt{\phi}_2=(\Phi(F_1)\circ \phi)^*(g')\circ Q_2(\Phi(F_1)\circ \phi)=Q_1(\Phi(F_1)\circ\phi)\circ g'$$
where the first equality is from the definition of $\iota$ in \cite[Construction 5.2]{fromunivwithPi} and the second from (\ref{2015.04.08.eq3}). On the other hand
$$\iota_1\circ \Phi(F_1^*(g))\circ \Phi(Q_2(F_1))\circ\wt{\phi}_2=\iota_1\circ\Phi(F_1^*(g)\circ Q_2(F_1))\circ \wt{\phi}_2=\iota_1\circ \Phi(Q_1(F_1)\circ g)\circ \wt{\phi}_2=$$
$$\iota_1\circ \Phi(Q_1(F_1))\circ \Phi(g)\circ \wt{\phi}_2=\iota_1\circ \Phi(Q_1(F_1))\circ \wt{\phi}_1\circ g'=Q_1(\Phi(F_1)\circ \phi)\circ g'$$
This proofs the first equation. For the second equation we have:
$$(\Phi(F_1)\circ \phi)^*(g')\circ \iota_2\circ \Phi(p_{X,F_1,2})=(\Phi(F_1)\circ \phi)^*(g')\circ p_{\Phi(X),\Phi(F_1)\circ\phi,2}=p_{\Phi(X),\Phi(F_1)\circ\phi,1}$$
and
$$\iota_1\circ \Phi(F_1^*(g))\circ \Phi(p_{X,F_1,2})=\iota_1\circ \Phi(F_1^*(g)\circ p_{X,F_1,2})=\iota_1\circ \Phi(p_{X,F_1,1})=p_{\Phi(X),\Phi(F_1)\circ \phi,1}$$
This finishes the proof of the lemma. 
\end{myproof}
\begin{lemma}
\llabel{2015.04.06.l7}
Under the previous assumptions and notations the squares
$$
\begin{CD}
\Phi(I_{p_2}(V)) @>\chi_2(V)>> I_{p_2'}(\Phi(V))\\
@V{\Phi(I^g(V))} VV @VV {I^{g'}(\Phi(V))} V\\
\Phi(I_{p_1}(V)) @>\chi_1(V)>> I_{p_1'}(\Phi(V))
\end{CD}
$$
commute.
\end{lemma}
\begin{myproof}
Let $X=I_{p_2}(V)$. We have
$$\chi_2(V)\circ {I^{g'}(\Phi(V))}=\eta'({\bf\Phi}_2^2(\eta^{-1}(Id_{X})))\circ I^{g'}(\Phi(V))$$
by definition of $\chi$ in \cite[Construction 5.6]{fromunivwithPi}. Then by Lemma \ref{2015.04.02.l4}  and Lemma \ref{2015.04.08.l1} we have:
$$\eta'({\bf\Phi}_2^2(\eta^{-1}(Id_{X})))\circ I^{g'}(\Phi(V))=\eta(D^{g'}(\Phi(X),\Phi(V))({\bf\Phi}_2^2(\eta^{-1}(Id_X))))=$$$$\eta({\bf\Phi}_1^2(D^g(X,V)(\eta^{-1}(Id_X))))$$
Then, again by Lemma \ref{2015.04.02.l4},   we have
$$\eta({\bf\Phi}_1^2(D^g(X,V)(\eta^{-1}(Id_X))))=\eta({\bf\Phi}_1^2(\eta^{-1}(Id_X\circ I^g(V)))=\eta({\bf\Phi}_1^2(\eta^{-1}(I^g(V))))$$
It remains to show that
$$\eta({\bf\Phi}_1^2(\eta^{-1}(I^g(V))))=\Phi(I^g(V))\circ \chi_1(V)$$
Let $a$ be any element of $Hom(I_{p_2}(V),I_{p_1}(V))$. Let us show that
$$\eta({\bf\Phi}_1^2(\eta^{-1}(a)))=\Phi(a)\circ \chi_1(V)$$
We have
$$\Phi(a)\circ \chi_1(V)=\Phi(a)\circ \eta'({\bf\Phi}^2_1(\eta^{-1}(Id_{I_{p_1}(V)})))=$$$$\eta'(D_{p_1'}(\Phi(a),\Phi(V))({\bf\Phi}^2_1(\eta^{-1}(Id_{I_{p_1}(V)})))$$
where the second equality holds because of naturality of $\eta$ in the first argument. Then 
$$\eta'(D_{p_1'}(\Phi(a),\Phi(V))({\bf\Phi}^2_1(\eta^{-1}(Id_{I_{p_1}(V)})))=\eta'({\bf\Phi}^2_1(D_{p_1}(a,V)(\eta^{-1}(I_{p_1}(V))))$$
by \cite[Lemma 5.4]{fromunivwithPi} and 
$$\eta'({\bf\Phi}^2_1(D_{p_1}(a,V)(\eta^{-1}(I_{p_1}(V))))=\eta'({\bf\Phi}^2_1(\eta^{-1}(a\circ Id_{p_1}(V))))=\eta'({\bf\Phi}_1^2(\eta^{-1}(a)))$$
again by naturality of $\eta$ in the first argument. This finishes the proof of Lemma \ref{2015.04.06.l7}.
\end{myproof}

Consider the morphisms
$$\zeta_i:\Phi(I_{p_i}(U))\sr I_{p_i}(U')$$
given by $\zeta_i=\chi_i(U)\circ I_{p_i}(\phi)$ and
$$\wt{\zeta}_i:\Phi(I_{p_i}(\wt{U}_1))\sr I_{p_i}(\wt{U}_1')$$
given by $\wt{\zeta}_i=\chi_i(\wt{U}_1)\circ I_{p_i}(\wt{\phi}_1)$. Note that 
$$\zeta_i=\xi_{(\Phi,\phi,\wt{\phi}_i)}$$
where $\xi$ are the morphisms introduced in \cite{fromunivwithPi} and $\wt{\zeta}_1=\xi_{(\Phi,\phi,\wt{\phi}_1)}$ but $\wt{\zeta}_2\ne \xi_{(\Phi,\phi,\wt{\phi}_2)}$.

\begin{theorem}
\llabel{2015.04.10.th3}
Under the previous assumptions and notations the morphisms $\zeta_1,\zeta_2,\wt{\zeta}_1,\wt{\zeta}_2$ form a morphism from the square
$$
\begin{CD}
\Phi(I_{p_2}(\wt{U}_1)) @>\Phi(I^g(\wt{U}_1))>> \Phi(I_{p_1}(\wt{U}_1))\\
@V\Phi(I_{p_2}(p_2)) VV @VV \Phi(I_{p_1}(p_2)) V\\
\Phi(I_{p_2}(U)) @>\Phi(I^g(U))>> \Phi(I_{p_1}(U))
\end{CD}
$$
to the square
$$
\begin{CD}
I_{p_2'}(\wt{U}_1') @>I^{g'}(\wt{U}_1')>> I_{p_1'}(\wt{U}_1')\\
@VI_{p_2'}(p_2') VV @VV I_{p_1'}(p_2') V\\
I_{p_2'}(U') @>I^{g'}(U')>> I_{p_1'}(U')
\end{CD}
$$
\end{theorem}
\begin{myproof}
We need to prove commutativity of the outer squares of the following four diagrams:
$$
\begin{CD}
\Phi(I_{p_2}(\wt{U}_1)) @>\chi_2(\wt{U}_1)>> I_{p_2'}(\Phi(\wt{U}_1)) @>I_{p_2'}(\wt{\phi}_1)>> I_{p_2'}(\wt{U}_1')\\
@V\Phi(I^g(\wt{U}_1))VV @VV I^{g'}(\Phi(\wt{U}_1)) V @VV I^{g'}(\wt{U}_1') V\\
\Phi(I_{p_1}(\wt{U}_1)) @>\chi_1(\wt{U}_1)>> I_{p_1'}(\Phi(\wt{U}_1)) @>I_{p_1'}(\wt{\phi}_1)>> I_{p_1'}(\wt{U}_1')
\end{CD}
$$
$$
\begin{CD}
\Phi(I_{p_2}(U_1)) @>\chi_2(U_1)>> I_{p_2'}(\Phi(U_1)) @>I_{p_2'}(\wt{\phi}_1)>> I_{p_2'}(U_1')\\
@V\Phi(I^g(U_1))VV @VV I^{g'}(\Phi(U_1)) V @VV I^{g'}(U_1') V\\
\Phi(I_{p_1}(U_1)) @>\chi_1(U_1)>> I_{p_1'}(\Phi(U_1)) @>I_{p_1'}(\wt{\phi}_1)>> I_{p_1'}(U_1')
\end{CD}
$$
$$
\begin{CD}
\Phi(I_{p_2}(\wt{U}_1)) @>\chi_2(\wt{U}_1)>> I_{p_2'}(\Phi(\wt{U}_1)) @>I_{p_2'}(\wt{\phi}_1)>> I_{p_2'}(\wt{U}_1')\\
@V\Phi(I_{p_2}(p_1)) VV @VV I_{p_2'}(\Phi(p_1)) V @VVI_{p_2'}(p_1') V\\
\Phi(I_{p_2}(U_1)) @>\chi_2(U_1)>> I_{p_2'}(\Phi(U_1)) @>I_{p_2'}(\wt{\phi}_1)>> I_{p_2'}(U_1')
\end{CD}
$$
$$
\begin{CD}
\Phi(I_{p_1}(\wt{U}_1)) @>\chi_1(\wt{U}_1)>> I_{p_1'}(\Phi(\wt{U}_1)) @>I_{p_1'}(\wt{\phi}_1)>> I_{p_2'}(\wt{U}_1')\\
@V\Phi(I_{p_1}(p_1)) VV @VV I_{p_1'}(\Phi(p_1)) V @VVI_{p_1'}(p_1') V\\
\Phi(I_{p_1}(U_1)) @>\chi_1(U_1)>> I_{p_1'}(\Phi(U_1)) @>I_{p_1'}(\wt{\phi}_1)>> I_{p_1'}(U_1')
\end{CD}
$$
The left squares in the first and the second diagram are commutative by Lemma \ref{2015.04.06.l7}.

The left squares in the third and the fourth diagram are commutative by \cite[Lemma 5.7]{fromunivwithPi}.

The right hand side squares in the first and second diagram commute by Lemma \ref{2015.04.10.l2}. 

The right hand side squares of the third and the fourth diagram commute because $I_{p'_i}$ are functorial and therefore take commutative squares to commutative squares.

Theorem is proved.
\end{myproof}

\subsection{Universe category functors compatible with J-structures}

Let us define now conditions on functors of universe categories that reflect the idea of compatibility with the J0- J1- and J2-structures on the universes. Recall that for universe categories $({\mathcal C},p,pt)$, $({\mathcal C}',p',pt')$ a functor of universe categories is a triple $(\Phi,\phi,\wt{\phi})$ where $\Phi:{\mathcal C}\sr {\mathcal C}'$ is a functor that takes the canonical squares to pull-back squares and $pt$ to a final object and $\phi:\Phi(U)\sr U'$, $\wt{\phi}:\Phi(\wt{U})\sr \wt{U}'$ are morphisms such that the square
\begin{eq}\llabel{2015.04.06.eq10}
\begin{CD}
\Phi(\wt{U}) @>\wt{\phi}>> \wt{U}'\\
@V\Phi(p) VV @Vp' VV\\
\Phi(U) @>\phi>> U
\end{CD}
\end{eq}
is a pull-back square. For any functor of universe categories and $X\in{\mathcal C}$, $F:X\sr U$ the morphism
$$\Phi(p_{X,F})*(\Phi(Q(F))\circ\wt{\phi}):\Phi((X;F))\sr (\Phi(X);\Phi(F)\circ\phi)$$
is an isomorphism and we will denote it $\Phi_{X,F}$. Let $\Phi\wt{U}p$ be the composition
$$\Phi((\wt{U};p)) \stackrel{\Phi_{\wt{U},p}}{\lr} (\Phi(\wt{U});\Phi(p)\circ \phi)=(\Phi(\wt{U});\wt{\phi}\circ p')\stackrel{Q'(\wt{\phi},p')}{\lr} (\wt{U}';p')$$
We also have another description of this morphism given by the following lemma.
\begin{lemma}
\llabel{2015.04.10.l5}
One has:
$$\Phi\wt{U}p=(\Phi(p_{\wt{U},p})\circ\wt{\phi})*(\Phi(Q(p))\circ\wt{\phi})$$
\end{lemma}
\begin{myproof}
One has
$$\Phi\wt{U}p\circ p'_{\wt{U}',p'}=\Phi_{\wt{U},p}\circ Q'(\wt{\phi},p')\circ p'_{\wt{U}',p'} =\Phi_{\wt{U},p}\circ p_{\Phi(\wt{U}),\wt{\phi}\circ p'}\circ\wt{\phi}=\Phi(p_{\wt{U},p})\circ \wt{\phi}$$
where the second equality is by definition of $Q'(-,-)$ and the third equality is by definition of $\Phi_{\wt{U},p}$. Then
$$\Phi\wt{U}p\circ Q'(p')=\Phi_{\wt{U},p}\circ Q'(\wt{\phi},p')\circ Q'(p')=\Phi_{\wt{U},p}\circ Q(\wt{\phi}\circ p')=$$$$\Phi_{\wt{U},p}\circ Q(\Phi(p)\circ \phi)=\Phi(Q(p))\circ \wt{\phi}$$
where again the second equality is by definition of $Q(-,-)$ and the fourth equality is by definition of $\Phi_{\wt{U},p}$. 
\end{myproof}
\begin{lemma}
\llabel{2015.04.10.l6}
For $s,s':Y\sr \wt{U}$ such that $s\circ p=s'\circ p$ one has
$$\Phi(s*s')\circ\Phi\wt{U}p=\Phi(s\circ \wt{\phi})*\Phi(s'\circ\wt{\phi})$$
in particular
$$\Phi(\Delta)\circ \Phi\wt{U}p = \wt{\phi}*\wt{\phi}$$
\end{lemma}
\begin{myproof}
Using Lemma \ref{2015.04.10.l5} we have
$$\Phi(s*s')\circ \Phi\wt{U}p\circ p'_{\wt{U}',p'}=\Phi(s*s')\circ\Phi(p_{\wt{U},p})\circ\wt{\phi}=s\circ \wt{\phi}$$
and
$$\Phi(s*s')\circ \Phi\wt{U}p\circ Q'(p')=\Phi(s*s')\circ \Phi(Q(p))\circ \wt{\phi}=s'\circ\wt{\phi}$$
The particular case of $\Delta$ follows from the fact that $\Delta=Id_{\wt{U}}*Id_{\wt{U}}$. 
\end{myproof}
\begin{lemma}
\llabel{2015.04.06.l5}
The square
$$
\begin{CD}
\Phi((\wt{U};p)) @>\Phi\wt{U}p>> (\wt{U}';p')\\
@VV \Phi(p_{\wt{U},p}) V @VV p'_{\wt{U}',p'}V\\
\Phi(\wt{U}) @>\wt{\phi}>> \wt{U}'
\end{CD}
$$
is a pull-back square.
\end{lemma}
\begin{myproof}
This square is equal to the composition of two squares
$$
\begin{CD}
\Phi((\wt{U};p)) @>\Phi_{\wt{U},p}>> (\Phi(\wt{U});\wt{\phi}\circ p') @>Q'(\wt{\phi},p')>> (\wt{U}';p')\\
@V\Phi(p_{\wt{U},p}) VV @VV p_{\Phi(\wt{U}),\wt{\phi}\circ p'}V @VVp'_{\wt{U}',p'}V\\
\Phi(\wt{U}) @= \Phi(\wt{U}) @>\wt{\phi}>> \wt{U}'
\end{CD}
$$
The right hand side square is a pull-back square (\ref{2015.04.06.l0.sq}). The left hand side square is a pull-back square as a commutative square whose sides are isomorphisms. We conclude that the composition of these two squares is a pull-back square.
\end{myproof}
\begin{definition}
\llabel{2015.04.06.def4}
Let $Eq$ be a J0-structure on $p$ and $Eq'$ a J0-structure on $p'$. A universe category functor $(\Phi,\phi,\wt{\phi})$ is said to be compatible with $Eq$ and $Eq'$ if the square
\begin{eq}\llabel{2015.04.06.eq6}
\begin{CD}
\Phi((\wt{U};p)) @>\Phi(Eq)>> \Phi(U)\\
@V\Phi\wt{U}p VV @VV\phi V\\
(\wt{U}';p') @>Eq'>> U'
\end{CD}
\end{eq}
commutes.
\end{definition}
Let $Eq$, $Eq'$ be as above. Let $(\Phi,\phi,\wt{\phi})$ be a universe functor compatible with $Eq$, and $Eq'$. Define a morphism
$$\wt{\phi}_{E}:\Phi(E\wt{U})\sr E\wt{U}'=((\wt{U}';p'),Eq')$$
as $(\Phi(p_{(\wt{U};p),Eq})\circ \Phi\wt{U}p)*(\Phi(Q(Eq))\circ \wt{\phi})$. 
\begin{lemma}
\llabel{2015.04.06.l4}
Let $Eq$, $Eq'$ be as above. Let $(\Phi,\phi,\wt{\phi})$ be a universe functor compatible with $Eq$, and $Eq'$. Then the square
$$
\begin{CD}
\Phi(E\wt{U}) @>\wt{\phi}_{E}>> E\wt{U}'\\
@V\Phi(p_{(\wt{U};p),Eq}) VV @VV p_{(\wt{U}';p'),Eq'}V\\
\Phi((\wt{U};p)) @>\Phi\wt{U}p>> (\wt{U}';p')
\end{CD}
$$
is a pull-back square.
\end{lemma} 
\begin{myproof}
Consider the diagram
\begin{eq}\llabel{2015.04.06.eq8}
\begin{CD}
\Phi(E\wt{U}) @>\wt{\phi}_{E}>> E\wt{U}' @>Q(Eq')>> \wt{U}'\\
@V\Phi(p_{(\wt{U};p),Eq}) VV @VV p_{(\wt{U}';p'),Eq'}V @VVp'V\\
\Phi((\wt{U};p)) @>\Phi\wt{U}p>> (\wt{U}';p') @>Eq'>> U'
\end{CD}
\end{eq}
The outer square of this diagram is equal to the outer square of the diagram
\begin{eq}\llabel{2015.04.06.eq7}
\begin{CD}
\Phi(E\wt{U}) @>\Phi(Q(Eq))>> \Phi(\wt{U}) @>\wt{\phi}>> \wt{U}'\\
@V\Phi(p_{(\wt{U};p),Eq}) VV @VV \Phi(p)V @VVp'V\\
\Phi((\wt{U};p)) @>\Phi(Eq)>> \Phi(U) @>\phi>> U'
\end{CD}
\end{eq}
where the equality of the lower horizontal arrows follows from the commutativity of the square (\ref{2015.04.06.eq6}). The left hand side square of this diagram is a pull-back square because $\Phi$ takes canonical squares to pull-back squares. The right hand side square is a pull-back square by definition of a functor of universe categories. Therefore the outer square is a pull-back square. The right hand side square of (\ref{2015.04.06.eq8}) is a canonical square and therefore a pull-back square. We conclude that the left hand square of (\ref{2015.04.06.eq8}) is a pull-back square.
\end{myproof}
\begin{lemma}
\llabel{2015.04.06.l6}
Let $Eq$, $Eq'$ be as above. Let $(\Phi,\phi,\wt{\phi})$ be a functor of universe categories compatible with $Eq$, and $Eq'$. Then the square
\begin{eq}\llabel{2015.04.06.eq9}
\begin{CD}
\Phi(E\wt{U}) @>\wt{\phi}_{E}>> E\wt{U}'\\
@V\Phi(pE\wt{U})VV @VV pE\wt{U}' V\\
\Phi(U) @>\phi>> U'
\end{CD}
\end{eq}
is a pull-back square.
\end{lemma} 
\begin{myproof}
It follows from the fact that the square (\ref{2015.04.06.eq9}) is equal to the vertical composition of the 
squares of Lemmas \ref{2015.04.06.l4} and \ref{2015.04.06.l5} with the square (\ref{2015.04.06.eq10}).
\end{myproof}

\begin{definition}
\llabel{2015.04.06.def5}
Let $Eq$, $Eq'$ be as above and let $\Omega$, $\Omega'$ be J1-structures over $Eq$ and $Eq'$ respectively. A universe category functor $(\Phi,\phi,\wt{\phi})$ is said to be compatible with $\Omega$ and $\Omega'$ if the square
$$
\begin{CD}
\Phi(\wt{U}) @>\Phi(\Omega)>> \Phi(\wt{U})\\
@V\wt{\phi} VV @VV\wt{\phi} V\\
\wt{U}' @>\Omega'>> \wt{U}'
\end{CD}
$$
commutes.
\end{definition}
\begin{lemma}
\llabel{2015.04.10.l7}
Let $Eq,\Omega$ and $Eq',\Omega'$ be as above and let ${\bf\Phi}$ be a universe category functor compatible with $Eq,Eq'$ and $\Omega,\Omega'$. Then the square
$$
\begin{CD}
\Phi(\wt{U}) @>\wt{\phi}>> \wt{U}'\\
@V\Phi(\omega)VV @VV\omega' V\\
\Phi(E\wt{U}) @>\wt{\phi}_{E}>> E\wt{U}'
\end{CD}
$$
commutes.
\end{lemma}
\begin{myproof}
Since $E\wt{U}'=((\wt{U}';p');Eq')$ it is sufficient to verify that the compositions of the two paths in the square with $p_{(\wt{U}';p'),Eq'}$ and $Q(Eq')$ coincide. We have:
$$\wt{\phi}\circ\omega'\circ Q(Eq')=\wt{\phi}\circ\Omega'$$
by definition of $\omega'$. On the other hand
$$\Phi(\omega)\circ \wt{\phi}_{E}\circ Q(Eq')=\Phi(\omega)\circ \Phi(Q(Eq))\circ \wt{\phi}=\Phi(\Omega)\circ\wt{\phi}$$
where the first equation holds by definition of $\wt{\phi}_{E}$. The proof follows now from the assumption that ${\bf\Phi}$ is compatible with $\Omega$ and $\Omega'$. 
\end{myproof}

To formulate the condition of compatibility of a universe functor with full J-structures on $\mathcal C$ and $\mathcal C'$ we will use Theorem \ref{2015.04.10.th3}.

Let ${\bf\Phi}=(\Phi,\phi,\wt{\phi})$ be a functor of universe categories. In view of Lemma \ref{2015.04.06.l6}, if $\Phi$ is compatible with $Eq$ and $Eq'$ then the triple ${\bf\Phi}_{E}:=(\Phi,\phi,\wt{\phi}_{E})$ is a functor of universe categories as well. If, in addition, ${\bf\Phi}$ is compatible with $\Omega$ and $\Omega'$ then, by Lemma \ref{2015.04.10.l7}, morphisms $\omega$ and $\omega'$ satisfy the conditions on morphisms $g$ and $g'$ of Section \ref{twouniv}. 

Let 
$$\xi_{\bf\Phi}:\Phi(I_p(U))\sr I_{p'}(U')$$
$$\wt{\xi}_{\bf\Phi}:\Phi(I_p(\wt{U}))\sr I_{p'}(\wt{U}')$$
denote the compositions $\chi_{\bf\Phi}(U)\circ I_{p'}(\phi)$ and $\chi_{\bf\Phi}(\wt{U})\circ I_{p'}(\wt{\phi})$ and let 
$$\zeta_{\bf\Phi}:\Phi(I_{pE\wt{U}}(U))\sr I_{pE\wt{U}'}(U')$$
$$\wt{\zeta}_{\bf\Phi}:\Phi(I_{pE\wt{U}}(\wt{U}))\sr I_{pE\wt{U}'}(\wt{U}')$$
be given by the compositions $\chi_{{\bf\Phi}_{E}}(U)\circ I_{pE\wt{U}'}(\phi)$ and $\chi_{{\bf\Phi}_{E}}(\wt{U})\circ I_{pE\wt{U}'}(\wt{\phi})$. Note that $\zeta_{\bf\Phi}=\xi_{{\bf\Phi}_{E}}$ but $\wt{\zeta}_{\bf\Phi}$ is different from $\wt{\xi}_{{\bf\Phi}_{E}}$ since the latter is equal to the composition $\chi_{{\bf\Phi}_{E}}(E\wt{U})\circ I_{pE\wt{U}'}(\wt{\phi}_{E})$. Applying Theorem \ref{2015.04.10.th3} in this context we get the following. 
\begin{theorem}
\llabel{2015.04.10.th1}
Let ${\bf\Phi}$ be a functor of universe categories compatible with the J1-structures $(Eq,\Omega)$ and $(Eq',\Omega')$ on $p$ and $p'$ respectively. Then the morphisms $\xi_{\bf\Phi}, \wt{\xi}_{\bf\Phi}, \zeta_{\bf\Phi}, \wt{\zeta}_{\bf\Phi}$ form a morphism from the square
$$
\begin{CD}
\Phi(I_{pE\wt{U}}(\wt{U})) @>\Phi(I^{\omega}(\wt{U}))>> \Phi(I_p(\wt{U}))\\
@V\Phi(I_{pE\wt{U}}(p))VV @VV\Phi(I_p(p))V\\
\Phi(I_{pE\wt{U}}(U)) @>\Phi(I^{\omega}(U))>> \Phi(I_p(U))
\end{CD}
$$
to the square
$$
\begin{CD}
I_{pE\wt{U}'}(\wt{U}') @>I^{\omega'}(\wt{U}')>> I_{p'}(\wt{U}')\\
@VI_{pE\wt{U}'}(p')VV @VVI_{p'}(p')V\\
I_{pE\wt{U}'}(U') @>I^{\omega}(U')>> I_{p'}(U')
\end{CD}
$$
\end{theorem}
Let $R_{\bf\Phi}$ denote the composition
$$\Phi((I_{pE\wt{U}}(U), I^{\omega}(U))\times_{I_p(U)} (I_p(\wt{U}), I_p(p)))\sr \Phi(I_{pE\wt{U}}(U), I^{\omega}(U))\times_{\Phi(I_p(U))}\Phi(I_p(\wt{U}), I_p(p))\sr$$
$$(I_{pE\wt{U}'}(U'),I^{\omega'}(U'))\times_{I_{p'}(U')}(I_{p'}(\wt{U}'),I_{p'}(p'))$$
where the second arrow is defined by $\xi_{\bf\Phi}, \wt{\xi}_{\bf\Phi}$ and  $\zeta_{\bf\Phi}$ in view of Theorem \ref{2015.04.10.th1}. 
\begin{definition}
\llabel{2015.04.06.def6}
Let $Eq$, $Eq'$, $\Omega$ and $\Omega'$ be as above. Let $Jp$ and $Jp'$ be J2-structures over $(Eq,\Omega)$ and $(Eq',\Omega')$ respectively.  A universe category functor $(\Phi,\phi,\wt{\phi})$ is said to be compatible with $Jp$ and $Jp'$ if it is compatible with $Eq$, $Eq'$ and $\Omega$, $\Omega'$ in the sense of Definitions \ref{2015.04.06.def4} and \ref{2015.04.06.def5} respectively and 
the square
$$
\begin{CD}
\Phi((I_{pE\wt{U}}(U), I^{\omega}(U))\times_{I_p(U)} (I_p(\wt{U}), I_p(p))) @>R_{\bf\Phi}>> (I_{pE\wt{U}'}(U'),I^{\omega'}(U'))\times_{I_{p'}(U')}(I_{p'}(\wt{U}'),I_{p'}(p'))\\
@V\Phi(Jp) VV @VVJp' V\\
\Phi(I_{pE\wt{U}}(\wt{U})) @>\wt{\zeta}_{\bf\Phi}>> I_{pE\wt{U}'}(\wt{U}')
\end{CD}
$$
commutes. 
\end{definition}

\subsection{Homomorphisms of C-systems compatible with J-structures}
\begin{definition}
\llabel{2015.04.06.def1}
Let $H:CC\sr CC'$ be a homomorphism of C-systems.
\begin{enumerate}
\item Let $IdT$, $IdT'$ are J0-structures on $CC$ and $CC'$ respectively.  Then $H$ is called a homomorphism of C-systems with J0-structures $(CC,IdT)\sr (CC,IdT')$ if for each $\Gamma\in Ob(CC)$ and $o,o'\in\wt{Ob}_1(\Gamma)$ such that $\partial(o)=\partial(o')$, one has 
$$H(IdT_{\Gamma}(o,o'))=IdT'_{H(\Gamma)}(H(o),H(o'))$$
(the right hand side of the equality makes sense because $H$ commutes with $\partial$).
\item Let $IdT$, $IdT'$ be as above and let $refl$, $refl'$ be J1-structures over $IdT$ and $IdT'$ respectively. A homomorphism of C-systems with J0-structures $H:(CC,IdT)\sr (CC',IdT')$ is called a homomorphism of C-systems with J1-structures 
$$(CC,IdT,refl)\sr (CC',IdT',refl')$$
if for all $\Gamma\in Ob(CC)$ and $o\in \wt{Ob}_1(\Gamma)$ one has
$$H(refl(o))=refl'(H(o))$$
\end{enumerate}
\end{definition}
For a C-system $CC$ with a J0-structure $IdT$ and a J1-structure $refl$ over $IdT$ define $Jdom(CC,IdT,refl)$ as the set of quadruples $(\Gamma,T,P,s0)$ where $\Gamma\in Ob$, $T\in Ob_1(\Gamma)$, $P\in Ob_1(IdxT(T))$ and $s0\in \wt{Ob}(rf_T^*(P))$. Equivalently we can say that $Jdom(CC,IdT,refl)$ is the subset in $Ob\times Ob\times Ob\times \wt{Ob}$ that consists of quadruples $(\Gamma,T,P,s0)$ where $ft(T)=\Gamma$, $ft(P)=IdxT(T)$ and $\partial(s0)=rf_T^*(P)$. Then a J2-structure is defined by a map $Jdom\sr \wt{Ob}$ with some properties.
\begin{lemma}
\llabel{2015.04.06.l3}
Let $H:CC\sr CC'$ be a homomorphism of C-systems. Let $\Gamma,X,Y\in Ob(CC)$, $m,n\in\nn$ and suppose that $ft^m(X)=ft^{n}(Y)=\Gamma$. Let $f:X\sr Y$ be a morphism over $\Gamma$ and let $F:\Gamma'\sr \Gamma$ be a morphism. Then 
$$H(F^*(f))=H(F)^*(H(f))$$
\end{lemma}
\begin{myproof}
This is easy to show from the defining properties of $F^*(f)$ and $H(F)^*(H(f))$.
\end{myproof}
\begin{lemma}
\llabel{2015.04.06.l2}
Let $IdT$, $IdT'$, $refl$ and $refl'$ be as in Definition \ref{2015.04.06.def1} and let 
$$H:(CC,IdT,refl)\sr (CC',IdT',refl')$$
be a homomorphism of C-systems with J1-structures. Then for all elements $(\Gamma,T,P,s0)$ of  $Jdom(IdT,refl)$ one has $(H(\Gamma),H(T),H(P),H(s0))\in Jdom(IdT',refl')$.
\end{lemma}
\begin{myproof}
We have $ft(H(T))=H(ft(T))=H(\Gamma)$ and $ft(H(P))=H(ft(P))=H(IdxT(T))$. We also have $\partial(H(s0))=H(\partial(s0))=H(rf_T^*(P))$. By Lemma \ref{2015.04.06.l3} we further have
$$H(rf_T^*(P))=H(rf_T)^*(H(P))$$
It remains to show that $H(IdxT(T))=IdxT'(H(T))$ and $H(rf_T)=rf_{H(T)}'$. This follows by a straightforward but lengthy computation from the defining equations (\ref{2015.04.06.eq1}) and (\ref{2015.04.02.eq1}).
\end{myproof}
\begin{definition}
\llabel{2015.04.06.def2}
Let $IdT$, $IdT'$, $refl$ and $refl'$ be as in Definition \ref{2015.04.06.def1} and let $J$, $J'$ be J2-structures over $(IdT,refl)$ and $(IdT',refl')$ respectively. A homomorphism of C-systems with J1-structures 
$$H:(CC,IdT,refl)\sr (CC',IdT',ref')$$
is called a homomorphism of C-systems with J-structures 
$$(CC,IdT,refl,J)\sr (CC',IdT',ref;',J')$$
if for all $\Gamma\in Ob(CC)$, $T\in Ob_1(\Gamma)$, $P\in Ob_1(IdxT(T))$ and $s0\in \wt{Ob}(rf_T^*(P))$ one has
$$H(J(\Gamma,T,P,s0))=J'(H(\Gamma),H(T),H(P),H(s0))$$
where the right hand side of the equation makes sense by Lemma \ref{2015.04.06.l2}.
\end{definition}

\subsection{Functoriality of the J-structures $(IdT_{Eq},refl_{\Omega},J_{Jp})$}
\label{2015.04.12.sec1}
Let us first remind that by \cite[Construction 3.3]{Cfromauniverse} any universe category functor ${\bf \Phi}=(\Phi,\phi,\wt{\phi})$ defines a homomorphism of C-systems
$$H:CC({\mathcal C},p)\sr CC({\mathcal C}',p')$$
To define $H$ on objects, one defines by induction on $n$, for all $\Gamma\in Ob_n(CC({\mathcal C},p))$,  pairs $(H(\Gamma),\psi_{\Gamma})$ where $H(\Gamma)\in Ob(CC({\mathcal C}',p'))$ and $\psi_{\Gamma}$ is a morphism 
$$\psi_{\Gamma}:int'(H(\Gamma))\sr \Phi(int(\Gamma))$$
as follows. For $n=0$ one has $H(())=()$ and $\psi_{()}:pt'\sr \Phi(pt)$ is the unique morphism to a final object $\Phi(pt)$. For $(\Gamma,F)\in Ob_{n+1}$ one has 
$$H((\Gamma,F))=(H(\Gamma),\psi_{\Gamma}\circ\Phi(F)\circ \phi)$$
and $\psi_{(\Gamma,F)}$ is the unique morphisms $int'(H(\Gamma,F))\sr \Phi(int(\Gamma,F))$ such that
$$\psi_{(\Gamma,F)}\circ \Phi(Q(F))\circ\wt{\phi}=Q'(\psi_{\Gamma}\circ\Phi(F)\circ\phi)$$
and
$$\psi_{(\Gamma,F)}\circ \Phi(p_{\Gamma,F})=p_{H((\Gamma,F))}\circ \psi_{\Gamma}$$
Observe that $\psi_{\Gamma}$ is automatically an isomorphism. The action of $H$ on morphisms is given, for $f:\Gamma\sr\Gamma'$, by
$$H(f)=\psi_{\Gamma}\circ\Phi(f)\circ\psi_{\Gamma'}^{-1}$$

\begin{lemma}
\llabel{2015.04.12.l1}
Let ${\bf\Phi}$ be a universe category functor as above that is compatible with the J0-structures $Eq$ and $Eq'$ on $p$ and $p'$ respectively. Then the homomorphism of C-systems $H=H({\bf\Phi})$ is a homomorphism of C-systems with J0-structures relative to $IdT_{Eq}$ and $IdT_{Eq'}$.
\end{lemma}
\begin{myproof}
Let $IdT=IdT_{Eq}$ and $IdT'=IdT_{Eq'}$. We need to check that for all $\Gamma\in Ob(CC({\mathcal C},p))$ and $o,o'\in \wt{Ob}_1(\Gamma)$ such that $\partial(o)=\partial(o')$ one has
$$H(IdT(o,o'))=IdT'(H(o),H(o'))$$
Since 
$$\partial(H(IdT(o,o'))=H(\Gamma)$$
$$\partial(IdT'(H(o),H(o')))=ft(\partial(H(o)))=H(ft(\partial(o)))=H(\Gamma)$$
this is equivalent to
$$u_1(H(IdT(o,o')))=u_1(IdT'(H(o),H(o')))$$
By \cite[Lemma 6.1(1)]{fromunivwithPi} we have
$$u_1(H(IdT(o,o')))=\psi_{\Gamma}\circ\Phi(u_1(IdT(o,o')))\circ\phi=\psi_{\Gamma}\circ\Phi((\wt{u}_1(o)*\wt{u}_1(o'))\circ Eq)\circ\phi=$$
$$\psi_{\Gamma}\circ\Phi(\wt{u}_1(o)*\wt{u}_1(o'))\circ \Phi(Eq)\circ\phi=\psi_{\Gamma}\circ\Phi(\wt{u}_1(o)*\wt{u}_1(o'))\circ \Phi\wt{U}p\circ Eq'$$
By Lemma \ref{2015.04.10.l6} we have 
$$\psi_{\Gamma}\circ\Phi(\wt{u}_1(o)*\wt{u}_1(o'))\circ \Phi\wt{U}p\circ Eq'=\psi_{\Gamma}\circ((\Phi(\wt{u}_1(o))\circ\wt{\phi})*(\Phi(\wt{u}_1(o'))\circ\wt{\phi}))\circ Eq'$$
and \cite[Lemma 6.1(2)]{fromunivwithPi}
$$\psi_{\Gamma}\circ((\Phi(\wt{u}_1(o))\circ\wt{\phi})*(\Phi(\wt{u}_1(o'))\circ\wt{\phi}))\circ Eq'=$$
$$((\psi_{\Gamma}\circ\Phi(\wt{u}_1(o))\circ\wt{\phi})*(\psi_{\Gamma}\circ\Phi(\wt{u}_1(o'))\circ\wt{\phi}))\circ Eq'=$$
$$(\wt{u}_1(H(o))*\wt{u}_1(H(o')))\circ Eq'=IdT'(H(o),H(o))$$
Lemma is proved.
\end{myproof}
\begin{lemma}
\llabel{2015.04.12.l2}
Let ${\bf\Phi}$ be a universe category functor as above that is compatible with the (J0,J1)-structures $(Eq,\Omega)$ and $(Eq',\Omega')$ on $p$ and $p'$ respectively. Then the homomorphism of C-systems $H=H({\bf\Phi})$ is a homomorphism of C-systems with (J0,J1)-structures relative to $(IdT_{Eq},refl_{\Omega})$ and $(IdT_{Eq'},refl_{\Omega'})$.
\end{lemma}
\begin{myproof}
Let $refl=refl_{\Omega}$ and $refl'=refl_{\Omega'}$. The compatibility condition is
$$\Phi(\Omega)\circ\wt{\phi}=\wt{\phi}\circ\Omega'$$
We need to check that for $\Gamma\in Ob(CC({\mathcal C},p))$ and $s\in \wt{Ob}_1(\Gamma)$ one has
$$H(refl(s))=refl'(H(s))$$
By \cite[Lemma 6.1(2)]{fromunivwithPi} we have
$$H(refl(s))=H(\wt{u}_1^{-1}(\wt{u}_1(s)\circ\Omega))=\wt{u}_1^{-1}(\psi_{\Gamma}\circ\Phi(\wt{u}_1(s)\circ\Omega)\circ\wt{\phi})=\wt{u}_1^{-1}(\psi_{\Gamma}\circ\Phi(\wt{u}_1(s))\circ\Phi(\Omega)\circ \wt{\phi})=$$
$$=\wt{u}_1^{-1}(\psi_{\Gamma}\circ\Phi(\wt{u}_1(s))\circ\wt{\phi}\circ\Omega')=\wt{u}_1^{-1}(\wt{u}_1(H(s))\circ\Omega')=refl'(H(s)).$$
\end{myproof}

To prove the functoriality of the full J-structures we will need some lemmas first.

Recall that in \cite{Csubsystems} we let $p_{\Gamma,n}:\Gamma\sr ft^n(\Gamma)$ denote the composition of $n$ canonical projections $p_{\Gamma}\circ \dots\circ p_{ft^{n-1}(\Gamma)}$. 
\begin{lemma}
\llabel{2015.05.10.l1}
Let ${\bf\Phi}$ be a universe category functor and $\Gamma\in Ob(CC({\mathcal C},p))$ be such that $l(\Gamma)\ge n$. Then the square
$$
\begin{CD}
int'(H(\Gamma)) @>\psi_{\Gamma}>> \Phi(int(\Gamma))\\
@Vp_{H(\Gamma),n} VV @VV\Phi(p_{\Gamma,n})V\\
int'(ft^n(\Gamma)) @>\psi_{ft^n(\Gamma)}>> \Phi(int(ft^n(\Gamma)))
\end{CD}
$$
commutes.
\end{lemma}
\begin{myproof}
It follows by simple induction from the defining relation
$$\psi_{\Gamma}\circ \Phi(p_{\Gamma})=p_{H(\Gamma)}\circ \psi_{ft(\Gamma)}$$
of $\psi_{\Gamma}$. 
\end{myproof}

\begin{lemma}
\llabel{2015.05.06.l3}
Let $Eq, Eq'$ be J0-structures on $({\mathcal C},p)$ and $({\mathcal C}',p')$ and $\bf\Phi$ be a universe category functor compatible with $Eq$ and $Eq'$. The for all $\Gamma\in Ob(CC({\mathcal C},p))$, $T\in Ob_1(\Gamma)$, $P\in Ob_1(IdxT(T))$ and $o\in \wt{Ob}(P)$ one has:
\begin{enumerate}
\item $(u'_{1,H(\Gamma)}(H(T)), u_{1,IdxT'(H(T))}'(H(P)))$ is a well defined element of $D_{pE\wt{U}'}(\Phi(int(\Gamma)),U')$ and
$$(u'_{1,H(\Gamma)}(H(T)), u_{1,IdxT'(H(T))}'(H(P)))=$$$$D_{pE\wt{U}'}(\psi_{\Gamma},\wt{U}')(D_{pE\wt{U}'}(int(H(\Gamma)),\phi)({\bf\Phi}_{E}^2(u_{1,\Gamma}(T),u_{1,Idx(T)}(P))))$$
\item $(u'_{1,H(\Gamma)}(H(T)), \wt{u}_{1,IdxT'(H(T))}'(H(o)))$ is a well defined element of $D_{pE\wt{U}'}(\Phi(int(\Gamma)),\wt{U}')$ and
$$(u'_{1,H(\Gamma)}(H(T)), \wt{u}_{1,IdxT'(H(T))}'(H(o)))=$$$$D_{pE\wt{U}'}(\psi_{\Gamma},\wt{U}')(D_{pE\wt{U}'}(int(H(\Gamma)),\wt{\phi})({\bf\Phi}_{E}^2(u_{1,\Gamma}(T),\wt{u}_{1,Idx(T)}(o))))$$
\end{enumerate}
\end{lemma}
\begin{remark}\rm
Since $u_2(T)=(u_1(ft(T)),u_1(T))$ and $\wt{u}_2(s)=(u_1(ft(\partial(s))),\wt{u}_1(s))$ , this lemma is very similar to \cite[Lemma 6.1(3,4)]{fromunivwithPi} but its proof is much more involved because of the interaction of the two different universe functors. 
\end{remark}
\begin{myproof}
We will only consider the second assertion. The proof of the first one is similar and simpler.  

To prove that the pair $(u'_{1,H(\Gamma)}(H(T)), \wt{u}_{1,IdxT'(H(T))}'(H(o)))$ is a well defined element of $D_{pE\wt{U}'}(\Phi(int(\Gamma)),\wt{U}')$ we need to show that $ft(\partial(H(o)))=IdxT'(H(T))$ and that 
the source of $\wt{u}_{1,IdxT'(H(T))}'(H(o)))$ equals to 
$(int(H(\Gamma)); u'_1(H(T)))_{E}$, i.e., that
$$int'(IdxT'(H(T)))=(int(H(\Gamma)); u'_1(H(T)))_{E}$$
The former is a corollary of our assumptions and Lemma \ref{2015.04.12.l1} and the latter is a corollary of 
\cite[Problem 3.3(1)]{fromunivwithPi} and the first equation of Lemma \ref{2015.03.27.l1}. 

Let $X=int(\Gamma)$, $F=u_{1,\Gamma}(T)$ and $\wt{G}=\wt{u}_{1,Idx(T)}(o)$. By definitions we have 
$$D_{pE\wt{U}'}(\psi_{\Gamma},\_)(D_{pE\wt{U}'}(\_,\wt{\phi})({\bf\Phi}_{E}^2(F,\wt{G})))=D_{pE\wt{U}'}(\psi_{\Gamma},\_)(D_{pE\wt{U}'}(\_,\wt{\phi})(\Phi(F)\circ \phi, \iota\circ \Phi(\wt{G})))=$$
$$D_{pE\wt{U}'}(\psi_{\Gamma},\_)(\Phi(F)\circ \phi, \iota\circ \Phi(\wt{G})\circ\wt{\phi})=(\psi_{\Gamma}\circ\Phi(F)\circ \phi, Q(\psi_{\Gamma},\Phi(F)\circ\phi)_{E'}\circ\iota\circ \Phi(\wt{G})\circ\wt{\phi})$$
where 
$$\iota:(\Phi(X);\Phi(F)\circ \phi)_{E'}\sr \Phi((X;F)_{E})$$
is the unique morphism such that
$$\iota\circ \Phi(p^E_{X,F})=p^{E'}_{\Phi(X),\Phi(F)\circ\phi}$$
$$\iota\circ \Phi(Q(F)_E)\circ \wt{\phi}_E=Q(\Phi(F)\circ\phi)_{E'}$$

On the other hand 
$$u_{1,H(\Gamma)}(H(T))=\psi_{\Gamma}\circ \Phi(u_{1,\Gamma}(T))\circ \phi$$
$$\wt{u}_{1,IdxT'(H(T))}(H(o))=\wt{u}_{1,H(IdxT(T))}(H(o))=\psi_{H(IdxT(T))}\circ \Phi(\wt{u}_{1,IdxT(T)}(o))\circ \wt{\phi}$$
by \cite[Lemma 6.1(1,2)]{fromunivwithPi}. Therefore, to prove the lemma it is sufficient to show that
$$\psi_{IdxT(T)}=Q(\psi_{\Gamma},\Phi(F)\circ\phi)_{E'}\circ\iota$$
Both sides are morphisms with the codomain 
$$\Phi(int(IdxT(T)))=\Phi((X;F)_E)$$
and since ${\bf\Phi}_E$ is a universe category functor it is sufficient to show that the compositions of the two sides with $\Phi(p^E_{X,F})$ and $\Phi(Q(F)_E)\circ \wt{\phi}_E$ are the same. Since
$$Q(\psi_{\Gamma},\Phi(F)\circ\phi)_{E'}\circ\iota\circ \Phi(p^E_{X,F})=Q(\psi_{\Gamma},\Phi(F)\circ\phi)_{E'}\circ p^{E'}_{\Phi(X),\Phi(F)\circ\phi}=$$
$$p^{E'}_{int'(H(\Gamma)),\psi_{\Gamma}\circ \Phi(F)\circ\phi}\circ\psi_{\Gamma}=p^{E'}_{int'(H(T))}\circ\psi_{\Gamma}$$
the first equation reduces to
\begin{eq}
\llabel{2015.05.10.eq1}
\psi_{IdxT(T)}\circ \Phi(p^E_{T})=p^{E'}_{int'(H(T))}\circ\psi_{\Gamma}
\end{eq}
and since
$$Q(\psi_{\Gamma},\Phi(F)\circ\phi)_{E'}\circ\iota\circ \Phi(Q(F)_E)\circ \wt{\phi}_E=Q(\psi_{\Gamma},\Phi(F)\circ\phi)_{E'}\circ Q(\Phi(F)\circ\phi)_{E'}=$$$$Q(\psi_{\Gamma}\circ \Phi(F)\circ \phi)_{E'}$$
the second equation reduces to
\begin{eq}
\llabel{2015.05.10.eq2.0}
\psi_{IdxT(T)}\circ \Phi(Q(F)_E)\circ \wt{\phi}_E=Q(\psi_{\Gamma}\circ \Phi(F)\circ \phi)_{E'}
\end{eq}
Equation (\ref{2015.05.10.eq1}) follows immediately from Lemma \ref{2015.05.10.l1} and the second equation of Lemma \ref{2015.03.27.l1}.

We have 
$$\wt{\phi}_E=(\Phi(p_{(\wt{U};p),Eq})\circ \Phi\wt{U}p)*(\Phi(Q(Eq))\circ\wt{\phi})$$
Therefore (\ref{2015.05.10.eq2.0}) is equivalent to two equations:
$$\psi_{IdxT(T)}\circ \Phi(Q(F)_E)\circ \Phi(p_{(\wt{U};p),Eq})\circ \Phi\wt{U}p=$$
\begin{eq}
\label{2015.05.10.eq2a}
Q(\psi_{\Gamma}\circ \Phi(F)\circ \phi)_{E'}\circ p_{(\wt{U}';p'),Eq'}
\end{eq}
and
$$\psi_{IdxT(T)}\circ \Phi(Q(F)_E)\circ \Phi(Q(Eq))\circ\wt{\phi}=$$
\begin{eq}
\label{2015.05.10.eq2b}
Q'(\psi_{\Gamma}\circ \Phi(F)\circ \phi)_{E'}\circ Q(Eq')
\end{eq}
The first equality we will have to decompose further into two using the fact that by Lemma \ref{2015.04.10.l5} we have 
$$\Phi\wt{U}p=(\Phi(p_{\wt{U},p})\circ\wt{\phi})*(\Phi(Q(p))\circ \wt{\phi})$$
Therefore (\ref{2015.05.10.eq2a}) is equivalent to two equations
$$\psi_{IdxT(T)}\circ \Phi(Q(F)_E)\circ \Phi(p_{(\wt{U};p),Eq})\circ \Phi(p_{\wt{U},p})\circ\wt{\phi}=$$
\begin{eq}
\llabel{2015.05.10.eq2aa}
Q(\psi_{\Gamma}\circ \Phi(F)\circ \phi)_{E'}\circ p_{(\wt{U}';p'),Eq'}\circ p_{\wt{U}',p'}
\end{eq}
and
$$\psi_{IdxT(T)}\circ \Phi(Q(F)_E)\circ \Phi(p_{(\wt{U};p),Eq})\circ \Phi(Q(p))\circ\wt{\phi}=$$
\begin{eq}
\llabel{2015.05.10.eq2ab}
Q(\psi_{\Gamma}\circ \Phi(F)\circ \phi)_{E'}\circ p_{(\wt{U}';p'),Eq'}\circ Q(p')
\end{eq}
To prove (\ref{2015.05.10.eq2b}) observe first two useful equalities
$$u_1(IdxT(T))=Q(Q(F),p)\circ Eq$$
$$Q(u_1(IdxT(T)))=Q(F)_{E}\circ Q(Eq)$$
where the first follows from the proof of Lemma \ref{2015.03.27.l1} and the second is the combination of the first with the third equality of the same lemma. 

Now we have:
$$\psi_{IdxT(T)}\circ \Phi(Q(F)_E)\circ \Phi(Q(Eq))\circ\wt{\phi}=\psi_{IdxT(T)}\circ \Phi(Q(F)_E\circ Q(Eq))\circ\wt{\phi}=$$
$$\psi_{IdxT(T)}\circ \Phi(Q(u_1(IdxT(T))))\circ \wt{\phi}=Q'(\psi_{IdxT(T)}\circ \Phi(u_1(IdxT(T)))\circ \phi)=Q(u_1(H(IdxT(T))))$$
and
$$Q(\psi_{\Gamma}\circ \Phi(F)\circ \phi)_{E'}\circ Q(Eq')=Q'(u_1(H(T)))_{E'}\circ Q(Eq')=Q'(u_1(H(IdxT(T))))$$
The equality (\ref{2015.05.10.eq2b}) is proved. 

To prove (\ref{2015.05.10.eq2aa}) observe two equalities:
$$Q(F)_E\circ p_{(\wt{U};p),Eq}=p_{IdxT(T)}\circ Q(Q(F),p)$$
$$Q(Q(F),p)\circ p_{\wt{U},p}=p_{ft(IdxT(T)}\circ Q(F)$$
The same equalities hold for $F'=\psi_{\Gamma}\circ\Phi(F)\circ \phi=u_1'(H(T))$ and the equation (\ref{2015.05.10.eq2aa}) becomes
$$\psi_{IdxT(T)}\circ \Phi(p_{IdxT(T)})\circ \Phi(p_{ft(IdxT(T))})\circ \Phi(Q(F))\circ \wt{\phi}=p_{IdxT'(H(T))}\circ p_{ft(IdxT'(H(T)))}\circ Q(F')$$
Using the defining equations for $\psi$ we rewrite the left hand side as
$$\psi_{IdxT(T)}\circ \Phi(p_{IdxT(T)})\circ \Phi(p_{ft(IdxT(T))})\circ \Phi(Q(F))\circ \wt{\phi}=p_{H(Idx(T))}\circ p_{ft(H(IdxT(T)))}\circ \psi_{T}\circ \Phi(Q(F))\circ \wt{\phi}$$
It remains to show that
$$Q(\psi_{\Gamma}\circ \Phi(F)\circ \phi)=\psi_{(\Gamma,F)}\circ \Phi(Q(F))\circ \wt{\phi}$$
which is the defining equation of $\psi_{(G,F)}$. 

To prove (\ref{2015.05.10.eq2ab}) let us rewrite the left hand side first
$$\psi_{IdxT(T)}\circ \Phi(Q(F)_E)\circ \Phi(p_{(\wt{U};p),Eq})\circ \Phi(Q(p))\circ \wt{\phi}=$$$$\psi_{IdxT(T)}\circ \Phi(p_{IdxT(T)}\circ Q(Q(F),p))\circ \Phi(Q(p))\circ\wt{\phi}=$$$$p_{H(IdxT(T))}\circ \psi_{ft(IdxT(T))}\circ \Phi(Q(Q(F)\circ p))\circ \wt{\phi}=$$$$p_{H(IdxT(T))}\circ \psi_{(\Gamma,F,Q(F)\circ p)}\circ \Phi(Q(Q(F)\circ p))\circ \wt{\phi}=p_{H(IdxT(T))}\circ Q(u_1(H(ft(IdxT(T)))))$$
Where the first equality holds in view of the upper square of Construction {2015.05.08.constr1}, the second one is one of the defining equalities of $\phi$ and the third one is from \cite[Lemma 3.2]{fromunivwithPi}. 

Let $F'=u_1(H(T))$. Rewriting the right hand side we get
$$Q(\psi_{\Gamma}\circ \Phi(F)\circ \phi)_{E'}\circ \circ p_{(\wt{U}';p'),Eq'}\circ Q(p')=Q(F')_{E'}\circ p_{(\wt{U}';p'),Eq'}\circ Q(p')=$$
$$p_{IdxT'(H(T)}\circ Q(Q(F'),p')\circ Q(p')=p_{IdxT'(H(T)}\circ Q(Q(F')\circ p')=$$
$$p_{IdxT'(H(T)}\circ Q(u_1(ft(IdxT(H(T)))))$$
Where the second equality is from the upper square of Construction {2015.05.08.constr1} in $\mathcal C'$, the third equality is from \cite[Lemma 3.2]{fromunivwithPi}, and the fourth from the middle square of Construction {2015.05.08.constr1} in $\mathcal C'$. 

Lemma \ref{2015.05.06.l3} is proved. 
\end{myproof}

\begin{lemma}
\llabel{2015.04.12.l3}
Let ${\bf\Phi}$ be a universe category functor as above that is compatible with the (J0,J1,J2)-structures $(Eq,\Omega,Jp)$ and $(Eq',\Omega',Jp)$ on $p$ and $p'$ respectively. Then the homomorphism of C-systems $H=H({\bf\Phi})$ is a homomorphism of C-systems with (J0,J1,J2)-structures relative to $(IdT_{Eq},refl_{\Omega},J_{Jp})$ and $(IdT_{Eq'},refl_{\Omega'},J_{Jp'})$.
\end{lemma}
\begin{myproof}
Let $IdT=IdT_{\Omega}$, $IdT'=IdT_{\Omega'}$, $refl=refl_{\Omega}$, $refl'=refl_{\Omega'}$, $J=J_{Jp}$ and $J'=J_{Jp'}$.

We need to verify that for all $\Gamma\in Ob(CC({\mathcal C},p))$, $T\in Ob_1(\Gamma)$, $P\in Ob_1(IdxT(T))$ and $s0\in \wt{Ob}(rf^*_T(P))$ one has
$$H(J(\Gamma,T,P,s0))=J'(H(\Gamma),H(T),H(P),H(s0))$$

The defining equation for $J'$ is
$$\eta_{pE\wt{U}'}(u'_1(H(T)),\wt{u}'_{1,{IdxT'(H(T)}}(J')))=\phi(H(\Gamma),H(T),H(P),H(s0))\circ Jp'$$
and to prove the lemma we need to show that $H(J)$ satisfies  this equation.

Using Lemma \ref{2015.05.06.l3} we have
$$\eta_{pE\wt{U}}(u'_1(H(T)),\wt{u}'_{1,{IdxT'(H(T)}}(H(J))))=$$$$\eta_{pE\wt{U}'}(D_{pE\wt{U}'}(\psi_{\Gamma},\_)(D_{pE\wt{U}'}(\_,\wt{\phi})(\Phi_E^2(u_1(T),\wt{u}_{1,IdxT(T)}(J)))))=$$
$$\psi_{\Gamma}\circ \eta_{pE\wt{U}'}(\Phi_E^2(u_1(T),\wt{u}_{1,IdxT(T)}(J)))\circ I_{pE\wt{U}'}(\wt{\phi})$$
By \cite[Lemma 5.8]{fromunivwithPi} and by the defining equation for $J$ we further have
$$\psi_{\Gamma}\circ \eta_{pE\wt{U}'}(\Phi_E^2(u_1(T),\wt{u}_{1,IdxT(T)}(J)))\circ I_{pE\wt{U}'}(\wt{\phi})=$$$$\psi_{\Gamma}\circ \Phi(\eta_{pE\wt{U}}(u_1(T),\wt{u}_{1,IdxT(T)}(J)))\circ \chi_{{\bf\Phi}_E}(\wt{U})\circ I_{pE\wt{U}'}(\wt{\phi})=$$
$$\psi_{\Gamma}\circ \Phi(\phi(\Gamma,T,P,s0)\circ Jp)\circ \chi_{{\bf\Phi}_E}(\wt{U})\circ I_{pE\wt{U}'}(\wt{\phi})=\psi_{\Gamma}\circ \Phi(\phi(\Gamma,T,P,s0)\circ Jp)\circ \wt{\zeta}_{\bf\Phi}$$

It remains to show that
$$\psi_{\Gamma}\circ \Phi(\phi(\Gamma,T,P,s0)\circ Jp)\circ \wt{\zeta}_{\bf\Phi}=\phi(H(\Gamma),H(T),H(P),H(s0))\circ Jp'$$
By the compatibility condition of Definition \ref{2015.04.06.def6} we see that it is sufficient to prove that 
$$\psi_{\Gamma}\circ \Phi(\phi(\Gamma,T,P,s0))\circ R_{\Phi} = \phi(H(\Gamma),H(T),H(P),H(s0))$$
Let 
$$pr_1:(I_{pE\wt{U}}(U),I^{\omega})\times_{I_{p}(U)}(I_{p}(\wt{U}),I_p(p))\sr I_{pE\wt{U}}(U)$$
$$pr_2:(I_{pE\wt{U}}(U),I^{\omega})\times_{I_{p}(U)}(I_{p}(\wt{U}),I_p(p))\sr I_{p}(\wt{U})$$
be the projections and let $pr_1'$, $pr_2'$ be their analogs in $\mathcal C'$. Then one has
$$R_{\Phi}\circ pr_1'=\Phi(pr_1)\circ \zeta_{\Phi}$$
$$R_{\Phi}\circ pr_2'=\Phi(pr_2)\circ \wt{\xi}_{\Phi}$$
On the other hand the defining relations of $\phi(\Gamma,T,P,s0)$ are
$$\phi(\Gamma,T,P,s0)\circ pr_1=\eta_{pE\wt{U}}(F,G)$$
$$\phi(\Gamma,T,P,s0)\circ pr_2=\eta_p(F,\wt{H})$$
where 
$$F=u_{1,\Gamma}(T)\spc\spc G=u_{1,IdxT(T)}(P)\spc\spc \wt{H}=\wt{u}_{1,T}(s0)$$
and similarly for $F',G'$ and $\wt{H}'$.

We need to prove
\begin{eq}
\llabel{2015.05.10.eq3a}
\psi_{\Gamma}\circ \Phi(\phi(\Gamma,T,P,s0))\circ R_{\Phi}\circ pr_1' = \phi(H(\Gamma),H(T),H(P),H(s0))\circ pr_1'
\end{eq}
and
\begin{eq}
\llabel{2015.05.10.eq3b}
\psi_{\Gamma}\circ \Phi(\phi(\Gamma,T,P,s0))\circ R_{\Phi}\circ pr_2' = \phi(H(\Gamma),H(T),H(P),H(s0))\circ pr_2'
\end{eq}
For (\ref{2015.05.10.eq3a}), rewriting the left hand side we get
$$\psi_{\Gamma}\circ \Phi(\phi(\Gamma,T,P,s0))\circ R_{\Phi}\circ pr_1' = \psi_{\Gamma}\circ \Phi(\phi(\Gamma,T,P,s0))\circ\Phi(pr_1)\circ \zeta_{\Phi}=$$
$$\psi_{\Gamma}\circ \Phi(\phi(\Gamma,T,P,s0)\circ pr_1)\circ \zeta_{\Phi}=\psi_{\Gamma}\circ \Phi(\eta_{pE\wt{U}}(F,G))\circ \zeta_{\Phi}$$
Continuing we get
$$\psi_{\Gamma}\circ \Phi(\eta_{pE\wt{U}}(F,G))\circ \zeta_{\Phi}=\psi_{\Gamma}\circ \Phi(\eta_{pE\wt{U}}(F,G))\circ \xi_{{\bf \Phi}_E}\circ I_{pE\wt{U}'}(\phi)=$$
$$\psi_{\Gamma}\circ\eta_{pE\wt{U}'}({\bf\Phi}^2_E(F,G))\circ I_{pE\wt{U}'}(\phi)$$
where the last equality holds by \cite[Lemma 5.8]{fromunivwithPi} applied to $X=int(\Gamma)$, $V=U$ and ${\bf\Phi}={\bf\Phi}_E$. Continuing further we get
$$\psi_{\Gamma}\circ\eta_{pE\wt{U}'}({\bf\Phi}^2_E(F,G))\circ I_{pE\wt{U}'}(\phi)=\eta_{pE\wt{U}'}(D_{pE\wt{U}}(\psi_{\Gamma},\_)(D_{pE\wt{U}}(\_,\phi)({\bf\Phi}^2(F,G))))=$$$$\eta_{pE\wt{U}'}((u_{1,H(\Gamma)}'(H(T)),u_{1,IdxT'(H(T))}(H(P))))$$
where the last equality holds by Lemma \ref{2015.05.06.l3}(1). 

Rewriting the right hand side we get
$$\phi(H(\Gamma),H(T),H(P),H(s0))\circ pr_1'=\eta_{pE\wt{U}'}(F',G')$$
where $F'=u_{1,H(\Gamma)}(H(T))$ and $G'=u_{1,IdxT'(H(\Gamma))}(H(P))$. 
This shows that the first equality holds.

For (\ref{2015.05.10.eq3b}), rewriting the left hand side we get
$$\psi_{\Gamma}\circ \Phi(\phi(\Gamma,T,P,s0))\circ R_{\Phi}\circ pr_2' = \psi_{\Gamma}\circ \Phi(\phi(\Gamma,T,P,s0))\circ \Phi(pr_2)\circ \wt{\xi}_{\Phi}=$$
$$\psi_{\Gamma}\circ \Phi(\phi(\Gamma,T,P,s0)\circ pr_2)\circ \wt{\xi}_{\Phi}=\psi_{\Gamma}\circ \Phi(\eta_p(F,\wt{H}))\circ \wt{\xi}_{\Phi}=\eta_{p'}(\wt{u}_{2,H(\Gamma)}'(H(s0)))$$
where the last equality holds by \cite[Lemma 6.2(2)]{fromunivwithPi} since $(F,\wt{H})=u_{2,\Gamma}(s0)$. 

Rewriting the right hand side we get
$$\phi(H(\Gamma),H(T),H(P),H(s0))\circ pr_2'=\eta_{p'}(F',\wt{H}')$$
which proves the second equality since $(F',\wt{H}')=\wt{u}_{2,H(\Gamma)}'(H(s0))$.
\end{myproof}

\def\cprime{$'$}


\end{document}